\documentclass[11pt, twoside]{article}
\usepackage{bcc-cambridge-class}
\usepackage{graphicx,tikz,txfonts}
\usetikzlibrary{calc}

\bcctitle{Transversals in Latin Squares}
\bccshorttitle{Transversals in Latin Squares}

\bccname{Richard Montgomery}
\bccshortname{Richard Montgomery}

\bccaddressa{Mathematics Institute\\ University of Warwick\\ Coventry\\ CV4 7AL\\ UK}
\bccemaila{richard.montgomery@warwick.ac.uk}
\begin{document}
%% BCC editors to set page numbers by uncommenting the line below and entering the correct number
%\setcounter{page}{<insert page number here>}

\makebcctitle

\begin{abstract}
A Latin square is an $n$ by $n$ grid filled with $n$ symbols so that each symbol appears exactly once in each row and each column. A transversal in a Latin square is a collection of cells which do not share any row, column, or symbol. This survey will focus on results from the last decade which have continued the long history of the study of transversals in Latin squares.
\end{abstract}

%%%%%%%%%%%%%%%%%%%%%%%%%%%%%%%%%%%%%%%%%%%%%%%%%%%%%%%%%%%%%%%%%%%%%%%%%%%%%%%%%%%%%%%%%%%%%%%%%%%%%%%%%%%%%%%%%%%%%%%%%%
%%%%%%%%%%%%%%%%%%%%%%%%%%%%%%%%%%%%%%%%%%%%%%%%%%%%%%%%%%%%%%%%%%%%%%%%%%%%%%%%%%%%%%%%%%%%%%%%%%%%%%%%%%%%%%%%%%%%%%%%%%
%%%%%%%%%%%%%%%%%%%%%%%%%%%%%%%%%%%%%%%%%%%%%%%%%%%%%%%%%%%%%%%%%%%%%%%%%%%%%%%%%%%%%%%%%%%%%%%%%%%%%%%%%%%%%%%%%%%%%%%%%%
%%%%%%%%%%%%%%%%%%%%%%%%%%%%%%%%%%%%%%%%%%%%%%%%%%%%%%%%%%%%%%%%%%%%%%%%%%%%%%%%%%%%%%%%%%%%%%%%%%%%%%%%%%%%%%%%%%%%%%%%%%

\section{Introduction}\label{sec:intro}
A \emph{Latin square of order $n$} is an $n$ by $n$ grid filled with $n$ symbols so that each symbol appears exactly once in each row and each column (see Figure~\ref{fig:LS}). Important examples of Latin squares include the multiplication tables of finite groups, where the rows and columns of the Latin square of order $n$ are both indexed by some group $G$ with order $n$, and each entry is the product of its row with its column. In part due to their connection with magic squares,
Latin squares have a long history dating back to early mathematics  (for more on this we recommend the excellent historical survey by Andersen~\cite{andersen2007history}). Euler initiated their more rigorous study in the 18th century, while, in modern mathematics, Latin squares have strong connections to design theory (as, for example, we will see in our discussion of Steiner triple systems) as well as links to 2-dimensional permutations, finite projective planes and error correcting codes (see~\cite{KPSY,wanless2011transversals} for more on this).

 A \emph{partial transversal} in a Latin square is a collection of cells which share no row, column or symbol, while a \emph{full transversal} is one which contains every symbol exactly once (sometimes referred to as simply a \emph{transversal}).
We will see that, when $n$ is even, there are numerous Latin squares which may not contain a full transversal; for now, we recall only the canonical such example by remarking that it is easy to show that the Latin square corresponding to the addition table for $\mathbb{Z}_{2m}$ does not contain a full transversal for any $m\in \mathbb{N}$ (where a group is abelian we will use addition rather than multiplication). Indeed, suppose to the contrary that it did contain a full transversal, $T$. Then, the sum of the entries in $T$ is the sum of the elements in $\mathbb{Z}_{2m}$, but also, by using the definition of the Latin square, the sum of the indices of the rows and the columns of the Latin square, and thus twice the sum of the elements of $\mathbb{Z}_{2m}$. That is, we have
\[
\sum_{v\in \mathbb{Z}_{2m}}v=\sum_{x\in T}\mathrm{entry}(x)=\sum_{x\in T}(\mathrm{row}(x)+\mathrm{column}(x))%=\sum_{v\in T}\mathrm{row}(v)+\sum_{v\in T}\mathrm{column}(v)
=2\sum_{v\in \mathbb{Z}_{2m}}v.
\]
As $\sum_{v\in \mathbb{Z}_{2m}}v=m$ in $\mathbb{Z}_{2m}$, this gives a contradiction, and thus the addition table for $\mathbb{Z}_{2m}$ contains no full transversal.

\begin{figure}
\begin{center}
\begin{tikzpicture}[define rgb/.code={\definecolor{mycolor}{rgb}{#1}}, rgb color/.style={define rgb={#1},mycolor},scale=1.2]

\def\wi{0.4cm}

\def\n{5}

\def\Latentry{{{A,K,Q,J}, {Q,J,A,K},{J,Q,K,A},{K,A,J,Q}}};

\foreach \x in {0,1,...,\n}
\foreach \y in {0,1,...,\n}
{
\coordinate (A\x\y) at ($\n*(0,\wi)+\x*(\wi,0)-\y*(0,\wi)$);
\draw (A\x\y) [thick] rectangle ($(A\x\y)+(\wi,\wi)$);
}
\foreach \x/\y in {0/0,1/1,3/2,5/3,2/4,4/5}
{
\draw (A\x\y) [thick,fill=black!40] rectangle ($(A\x\y)+(\wi,\wi)$);
}

\foreach \x/\y in {0/0,4/1,5/2,2/3,3/4,1/5}
{
\draw ($(A\x\y)+0.5*(\wi,\wi)$) node {1};
}
\foreach \x/\y in {5/0,1/1,2/2,0/3,4/4,3/5}
{
\draw ($(A\x\y)+0.5*(\wi,\wi)$) node {2};
}
\foreach \x/\y in {4/0,3/1,1/2,5/3,0/4,2/5}
{
\draw ($(A\x\y)+0.5*(\wi,\wi)$) node {3};
}
\foreach \x/\y in {1/0,2/1,3/2,4/3,5/4,0/5}
{
\draw ($(A\x\y)+0.5*(\wi,\wi)$) node {4};
}
\foreach \x/\y in {3/0,0/1,4/2,1/3,2/4,5/5}
{
\draw ($(A\x\y)+0.5*(\wi,\wi)$) node {5};
}
\foreach \x/\y in {2/0,5/1,0/2,3/3,1/4,4/5}
{
\draw ($(A\x\y)+0.5*(\wi,\wi)$) node {6};
}

\end{tikzpicture}\hspace{2cm}\begin{tikzpicture}[define rgb/.code={\definecolor{mycolor}{rgb}{#1}}, rgb color/.style={define rgb={#1},mycolor},scale=1.2]

\def\wi{0.4cm}

\def\n{5}

\def\Latentry{{{A,K,Q,J}, {Q,J,A,K},{J,Q,K,A},{K,A,J,Q}}};

\foreach \x in {0,1,...,\n}
\foreach \y in {0,1,...,\n}
{
\coordinate (A\x\y) at ($\n*(0,\wi)+\x*(\wi,0)-\y*(0,\wi)$);
\draw (A\x\y) [thick] rectangle ($(A\x\y)+(\wi,\wi)$);
}
\foreach \x/\y in {2/1,5/1,2/4,5/4}
{
\draw (A\x\y) [thick,fill=black!25] rectangle ($(A\x\y)+(\wi,\wi)$);
}

\foreach \x/\y in {0/0,5/1,4/2,3/3,2/4,1/5}
{
\draw ($(A\x\y)+0.5*(\wi,\wi)$) node {0};
}
\foreach \x/\y in {1/0,0/1,5/2,4/3,3/4,2/5}
{
\draw ($(A\x\y)+0.5*(\wi,\wi)$) node {1};
}
\foreach \x/\y in {2/0,1/1,0/2,5/3,4/4,3/5}
{
\draw ($(A\x\y)+0.5*(\wi,\wi)$) node {2};
}
\foreach \x/\y in {3/0,2/1,1/2,0/3,5/4,4/5}
{
\draw ($(A\x\y)+0.5*(\wi,\wi)$) node {3};
}
\foreach \x/\y in {4/0,3/1,2/2,1/3,0/4,5/5}
{
\draw ($(A\x\y)+0.5*(\wi,\wi)$) node {4};
}
\foreach \x/\y in {5/0,4/1,3/2,2/3,1/4,0/5}
{
\draw ($(A\x\y)+0.5*(\wi,\wi)$) node {5};
}

\end{tikzpicture}
\end{center}
\caption{A Latin square of order 6 with a full transversal highlighted, and one with no full transversal (the addition table for $\mathbb{Z}_6$). For later illustration, an \emph{intercalate}, or 2 by 2 Latin subsquare, is highlighted in the second Latin square.}\label{fig:LS}
\end{figure}
These examples already lead to two natural questions: Which Latin squares contain a full transversal? If a Latin square has no full transversal, then what is the size of a largest partial transversal? These are interesting and very challenging questions, and the progress on them will be the focus of this survey.

In part due to the historical context, we will start by considering Latin squares with a much stronger property than simply the existence of a full transversal, those Latin squares which can be decomposed entirely into disjoint full transversals (see Section~\ref{sec:decomp}), while also discussing random Latin squares. Recently, Latin squares have been often studied using an equivalent formulation in edge-coloured complete bipartite graphs, and in Section~\ref{sec:examples} we will recall this formulation before using it to describe a large class of Latin squares with no full transversal. We will then discuss (in Section~\ref{sec:full}) Latin squares which are known to have a full transversal, before considering (in Section~\ref{sec:partial}) large partial transversals in any Latin square and, in particular, the well-known Ryser-Brualdi-Stein conjecture. Finally, in Section~\ref{sec:related}, we will consider other problems related to the study of transversals in Latin squares.

Before continuing any further, we highlight two related surveys which, like the present survey, were written in connection with editions of the British Combinatorial Conference. Firstly, the 2011 survey by Wanless~\cite{wanless2011transversals} also on `Transversals in Latin squares';
while many of the questions and results covered by~\cite{wanless2011transversals} are also covered here with their subsequent developments, other questions are also considered or are covered in more detail. Secondly, we discuss here the links to the study of rainbow subgraphs in edge-coloured graphs most closely related to transversals in Latin squares, but many further problems and more detailed discussion can be found in the 2022 survey by Pokrovskiy~\cite{Alexeysurvey} on `Rainbow subgraphs and their applications'.

%%%%%%%%%%%%%%%%%%%%%%%%%%%%%%%%%%%%%%%%%%%%%%%%%%%%%%%%%%%%%%%%%%%%%%%%%%%%%%%%%%%%%%%%%%%%%%%%%%%%%%%%%%%%%%%%%%%%%%%%%%
%%%%%%%%%%%%%%%%%%%%%%%%%%%%%%%%%%%%%%%%%%%%%%%%%%%%%%%%%%%%%%%%%%%%%%%%%%%%%%%%%%%%%%%%%%%%%%%%%%%%%%%%%%%%%%%%%%%%%%%%%%
%%%%%%%%%%%%%%%%%%%%%%%%%%%%%%%%%%%%%%%%%%%%%%%%%%%%%%%%%%%%%%%%%%%%%%%%%%%%%%%%%%%%%%%%%%%%%%%%%%%%%%%%%%%%%%%%%%%%%%%%%%
%%%%%%%%%%%%%%%%%%%%%%%%%%%%%%%%%%%%%%%%%%%%%%%%%%%%%%%%%%%%%%%%%%%%%%%%%%%%%%%%%%%%%%%%%%%%%%%%%%%%%%%%%%%%%%%%%%%%%%%%%%

\section{Decompositions into full transversals and random Latin squares}\label{sec:decomp}
%The consideration of decompositions of Latin squares into full transversals dates back at least to a recreational mathematics problem popular from at least the 1720's when it appeared in the compendium `R\'ecr\'eations math\'ematiques et physiques' by Jacques Ozanam~\cite{ozanam1723recreations}.
A Latin square decomposed into full transversals appears in a recreational mathematics problem which was popular from at least the 1720s when it appeared in the compendium `R\'ecr\'eations math\'ematiques et physiques' by Jacques Ozanam~\cite{ozanam1723recreations}.
Taking a  standard deck of cards, can the aces, kings, queens and jacks be arranged into a 4 by 4 grid so that each column and each row contains an ace, king, queen and jack which are all from different suits?
The answer is yes, and indeed a solution is given in Figure~\ref{fig:soln}, displayed overleaf so as not to disappoint any enthusiastic readers. In such a solution, forgetting the suit of each card will give a Latin square of order 4, in which, furthermore, each suit marks out a full transversal. That is, it gives a Latin square of order $4$ which can be decomposed into full transversals.

Euler~\cite{euler} took this problem further around 1780, for Latin squares of order 6, in the form of his famous `36 officers problem'. In this problem, there are 36 officers of 6 different ranks from 6 different regiments, with one officer from each rank from each regiment, who are to stand in a 6 by 6 grid so that each row and column contains officers of different ranks from different regiments. To notate this, Euler used Latin letters to denote the regiments and Greek letters to denote the ranks; the Latin letters used thus form a Latin square, giving rise to the terminology `Latin square'. While each rank (or Greek letter) would also mark out a full transversal and thus give a decomposition of a Latin square of order 6 into full transversals, the grid of ranks (or Greek letters) would also form a Latin square. We would say that the two Latin squares at play here are \emph{orthogonal}: two Latin squares of order $n$ are said to be \emph{orthogonal} (sometimes \emph{mutually orthogonal}) if the pairs of entries which appear in corresponding cells are distinct (so that all possible $n^2$ pairs appear). Finding two orthogonal Latin squares of order $n$ is equivalent to finding one Latin square of order $n$ that can be decomposed into full transversals.

Euler wrote that he could find no solution to his 36 officers problem, though could not prove none exists, and conjectured that, for each $n\equiv 2\;\mathrm{mod}\; 4$, there is no Latin square of order $n$ which can be decomposed into disjoint full transversals~\cite{euler}. This can immediately be seen to be true for $n=2$, but it was not until 1900 that Tarry~\cite{tarry1900probleme} confirmed this for $n=6$, thus showing that Euler's 36 officers problem has no solution. More generally, however, Euler's conjecture is false! This was first shown by Bose and Shrikhande~\cite{BS} in 1959 who constructed counterexamples for $n=22$ and $n=50$, before extending this with Parker~\cite{BPS} to show that Euler's conjecture is false for any $n\equiv 2\;\mathrm{mod}\; 4$ with $n\geq 10$.

Given this, an interesting question is to ask how common counterexamples are, and, in particular, is a random Latin square of order $n$ likely to be a counterexample to Euler's conjecture? For each $n\in \mathbb{N}$, let $\mathcal{L}(n)$ be the set of Latin squares of order $n$ which use the symbols in $[n]=\{1,2,\ldots,n\}$, and let $L_n\in \mathcal{L}(n)$ be a Latin square chosen from $\mathcal{L}(n)$ uniformly at random.
In 1990, van Rees~\cite{vanrees1990subsquares} conjectured that a vanishingly small proportion of Latin squares have a decomposition into full transversals, or, equivalently, that the random Latin square $L_n$ has no decomposition into full transversals with high probability. Wanless and Webb~\cite{wanless2006existence}, however, observed in 2006 that numerical observations suggest the contrary, and, considering the recent advances in techniques and the results shown, the following conjecture now appears likely.

\begin{conj}\label{conj:random}
The random Latin square $L_n$ can be decomposed into full transversals with high probability.% $1-o(1)$.
\end{conj}

Kwan~\cite{kwan2020almost} showed in 2020 that the random Latin square $L_n$ at least has a full transversal with high probability. Here, the rigidity of Latin squares make it difficult to study almost any non-trivial properties of a uniformly random Latin square. Kwan~\cite{kwan2020almost} approached the random Latin square $L_n$ by approximating it via a modified random triangle removal process, adapting breakthrough methods of Keevash~\cite{keevash2014existence,keevash2018counting} on the existence of designs and the analysis of the random triangle removal process carried out by Bohman, Frieze and Lubetzky~\cite{bohman2015random}. To find a full transversal, Kwan~\cite{kwan2020almost} used the absorbing method, as codified by R\"odl, Ruci\'nski and Szemer\'edi~\cite{rodl2006dirac} in 2006 (following earlier related work by Erd\H{o}s, Gy\'arf\'as and Pyber~\cite{EGP} and Krivelevich~\cite{krivtri}). As used in many different settings since, the aim here is find a small partial transversal (functioning as an absorption structure) with some flexibility before extending it to an almost-full transversal. Finally, the flexibility of the absorption structure is used to make a relatively small adjustment to the overall structure to extend this to a full transversal, where each column, row and symbol must appear exactly once.

In this setting, it is typically not hard to extend some small partial transversal into an almost-full transversal in any Latin square. For some potential applications this can moreover be done so that the unused columns, rows and symbols look like a randomly chosen set, for example using the semi-random method (also known as the R\"odl nibble) as introduced by R\"odl~\cite{rodlandhisnibble}. This method was used, via a development by Pippenger and Spencer~\cite{pippenger1989asymptotic}, in random Latin squares by Ferber and Kwan~\cite{ferbkwan} (for the work discussed below), while it was implemented for generalised Latin squares (see Section~\ref{sec:related}) by Montgomery, Pokrovskiy and Sudakov~\cite{montgomery2018decompositions} who showed that any generalised Latin square contains $(1-o(1))n$ disjoint partial transversals of order $(1-o(1))n$. (For more on nibble methods, see also the recent survey by Kang, Kelly, K{\"u}hn, Methuku and Osthus~\cite{kang2021graph}.)

\begin{figure}
\begin{center}
\includegraphics[scale = 1.5]{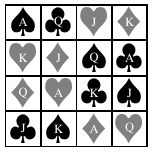}
\end{center}

\vspace{-0.2cm}

\caption{A Latin square of order 4 which can be decomposed into disjoint full transversals, as indicated by the card suits.}\label{fig:soln}
\end{figure}

The challenge then in using an absorption approach to find a full transversal in a random Latin square is to show that an appropriate absorption structure likely exists. To do this, Kwan~\cite{kwan2020almost} built up an absorption structure from much smaller flexible substructures. This uses distributive absorption, a specific style of absorption introduced by Montgomery~\cite{montgomery2018spanning} which has been applied in a range of settings to find large rigid substructures. More detail on distributive absorption is given in Section~\ref{sec:partial}.

%2
An alternative approach to studying full transversals in the random Latin square $L_n$ was introduced by Gould and Kelly~\cite{gould2023hamilton}, who showed moreover that the random Latin square $L_n$ has a \emph{Hamilton} transversal with high probability. Here, a \emph{Hamilton transversal} is a full transversal where the natural permutation of the indexing set of the rows/columns corresponding to the full transversal is a cyclic permutation. (For example, the full transversal highlighted in Figure~\ref{fig:LS} is not a Hamilton transversal as the corresponding natural permutation is $(1)(2)(3465)$.)
For their result, Gould and Kelly built on their work with K\"uhn and Osthus~\cite{gould2022almost} (as described briefly in Section~\ref{sec:related:cycle})
by using `switching methods' to study random Latin squares. Instead of working with random Latin squares, the idea is to consider
random \emph{Latin rectangles}, where a $k$ by $n$ Latin rectangle is a $k$ by $n$ array filled using $n$ symbols, so that each symbol is used at most one in each row  and column (see Figure~\ref{fig:rect}). As Latin rectangles are less rigid than Latin squares, small alterations (known as `switchings') can be made to a Latin rectangle in many different ways to reach another Latin rectangle. Gould and Kelly use the analysis of switching arguments to show that a uniformly random $k$ by $n$ Latin rectangle contains their desired small structures (for some $k=\Theta(n)$). These results can then be moved into a random Latin square by comparing its first $k$ rows to a uniformly random $k$ by $n$ Latin rectangle, an approach pioneered by McKay and Wanless~\cite{mckay1999most} which uses estimates on the permanent (more specifically an upper bound of Br\`egman~\cite{bregman1973some} and a lower bound due to Egorychev~\cite{egorychev1981solution} and Falikman~\cite{falikman1981proof}) to bound the number of extensions of a Latin rectangle to a Latin square.

\begin{figure}
\begin{center}
\begin{tikzpicture}[define rgb/.code={\definecolor{mycolor}{rgb}{#1}}, rgb color/.style={define rgb={#1},mycolor},scale=1.2]

\def\wi{0.4cm}

\def\n{5}

\def\Latentry{{{A,K,Q,J}, {Q,J,A,K},{J,Q,K,A},{K,A,J,Q}}};

\foreach \x in {0,1,...,\n}
\foreach \y in {0,1,2}
{
\coordinate (A\x\y) at ($\n*(0,\wi)+\x*(\wi,0)-\y*(0,\wi)$);
\draw (A\x\y) [thick] rectangle ($(A\x\y)+(\wi,\wi)$);
}
\foreach \x/\y in {0/0,1/1,0/1,1/0}
{
\draw (A\x\y) [thick,fill=black!25] rectangle ($(A\x\y)+(\wi,\wi)$);
}

\foreach \x/\y in {0/0,4/1,5/2}
{
\draw ($(A\x\y)+0.5*(\wi,\wi)$) node {1};
}
\foreach \x/\y in {5/0,1/1,2/2}
{
\draw ($(A\x\y)+0.5*(\wi,\wi)$) node {2};
}
\foreach \x/\y in {4/0,3/1,1/2}
{
\draw ($(A\x\y)+0.5*(\wi,\wi)$) node {3};
}
\foreach \x/\y in {1/0,2/1,3/2}
{
\draw ($(A\x\y)+0.5*(\wi,\wi)$) node {4};
}
\foreach \x/\y in {3/0,0/1,4/2}
{
\draw ($(A\x\y)+0.5*(\wi,\wi)$) node {5};
}
\foreach \x/\y in {2/0,5/1,0/2}
{
\draw ($(A\x\y)+0.5*(\wi,\wi)$) node {6};
}

\end{tikzpicture}\hspace{2cm}\begin{tikzpicture}[define rgb/.code={\definecolor{mycolor}{rgb}{#1}}, rgb color/.style={define rgb={#1},mycolor},scale=1.2]

\def\wi{0.4cm}

\def\n{5}

\def\Latentry{{{A,K,Q,J}, {Q,J,A,K},{J,Q,K,A},{K,A,J,Q}}};

\foreach \x in {0,1,...,\n}
\foreach \y in {0,1,2}
{
\coordinate (A\x\y) at ($\n*(0,\wi)+\x*(\wi,0)-\y*(0,\wi)$);
\draw (A\x\y) [thick] rectangle ($(A\x\y)+(\wi,\wi)$);
}
\foreach \x/\y in {0/0,1/1,0/1,1/0}
{
\draw (A\x\y) [thick,fill=black!25] rectangle ($(A\x\y)+(\wi,\wi)$);
}

\foreach \x/\y in {1/0,4/1,5/2}
{
\draw ($(A\x\y)+0.5*(\wi,\wi)$) node {1};
}
\foreach \x/\y in {5/0,0/1,2/2}
{
\draw ($(A\x\y)+0.5*(\wi,\wi)$) node {2};
}
\foreach \x/\y in {4/0,3/1,1/2}
{
\draw ($(A\x\y)+0.5*(\wi,\wi)$) node {3};
}
\foreach \x/\y in {0/0,2/1,3/2}
{
\draw ($(A\x\y)+0.5*(\wi,\wi)$) node {4};
}
\foreach \x/\y in {3/0,1/1,4/2}
{
\draw ($(A\x\y)+0.5*(\wi,\wi)$) node {5};
}
\foreach \x/\y in {2/0,5/1,0/2}
{
\draw ($(A\x\y)+0.5*(\wi,\wi)$) node {6};
}

\end{tikzpicture}
\end{center}
\caption{A 3 by 6 Latin rectangle on the left, where the comparative lack of rigidity allows many different small alterations to made while maintaining the Latin rectangles properties, with one such operation carried out on the highlighted cells to create the Latin rectangle on the right. Note that the same operation applied to the first Latin square in Figure~\ref{fig:LS} does not result in a Latin square.}\label{fig:rect}
\end{figure}

%3
Very recently, Eberhard, Manners and Mrazovi\'c~\cite{eberhard2023transversals} introduced yet another method for showing a random Latin square has a full transversal with high probability, giving a very different proof by using tools from analytic number theory.  In particular, their methods are extremely powerful for counting the number of transversals that can be expected in a typical random Latin square. In~\cite{kwan2020almost}, Kwan had shown that  $L_n$ has $\left((1-o(1))\frac{n}{e^2}\right)^n$ full transversals with high probability (and a similar counting result for Hamilton transversals was given by Gould and Kelly~\cite{gould2023hamilton}). Here, the upper bound on the number of full transversals is due to Taranenko~\cite{taranenko2015multidimensional} (with a simpler proof subsequently given by Glebov and Luria~\cite{glebov2016maximum}) and holds for any Latin square of order $n$.
Eberhard, Manners and Mrazovi\'c~\cite{eberhard2023transversals} gave the remarkably strong bound that a random Latin square has
$\big(e^{-1/2}+o(1)\big)\frac{(n!)^2}{n^n}$ full transversals with high probability. (See Section~\ref{sec:full} for a further discussion of their methods and results.)

%4
So far, the results mentioned sought one full transversal in a random Latin square (though finding many such single full transversals), some way from proving Conjecture~\ref{conj:random}!
Further work in this direction has taken place in the closely related setting of perfect matchings in random Steiner triple systems. This is described in Section~\ref{sec:related:STS}, but the results shown are likely to follow through with appropriate modification in Latin squares. That is, methods by Ferber and Kwan~\cite{ferbkwan} on Steiner triple systems are likely with modification to show that the random Latin square $L_n$ contains $(1-o(1))n$ disjoint full transversals with high probability (as noted in~\cite{ferbkwan}). These methods build on those by Kwan in \cite{kwan2020almost} by introducing sparse regularity techniques (as well as using a random partitioning argument of Ferber, Kronenberg and Long~\cite{ferber2017packing}). More specifically, they require a generalisation of the sparse regularity lemma of Kohayakawa and R\"odl~\cite{kohayakawa2003szemeredi} to hypergraphs, to show a sparse version of a `weak' hypergraph regularity lemma of Kohayakawa, Nagle, R\"odl and Schacht~\cite{kohayakawa2010weak}, which is then used in conjunction with a generalistion to linear hypergraphs of the resolution of the K\L R conjecture by Conlon, Gowers, Samotij and Schacht~\cite{conlon2014klr}.

The results mentioned here on random Latin squares have all been shown within the last decade, with most of them more recent still. These new techniques open up the prospect of a proof of Conjecture~\ref{conj:random} in the near future, though additional new ideas are no doubt needed. As highlighted by Ferber and Kwan~\cite{ferbkwan}, it is natural to take an absorption approach towards Conjecture~\ref{conj:random}. Most implementations of the absorption method in the literature correspond to finding a full transversal, but some more elaborate implementations have been used to decompose structures entirely. For example, Barber, Kühn, Lo and Osthus~\cite{barber2016edge} introduced iterative absorption techniques to decompose dense graphs into small fixed subgraphs and Glock, K{\"u}hn, Montgomery and Osthus~\cite{glock2021decompositions} decomposed large complete properly coloured graphs into rainbow spanning trees (see Section~\ref{sec:examples} for an exposition of the links of the current problem to rainbow subgraphs).

Finally, in this section, let us mention recent techniques to study the distribution of small structures within random Latin squares.
The results mentioned above using distributive absorption (i.e., \cite{ferbkwan,gould2023hamilton,kwan2020almost}) built their absorption structure from small structures found in Latin squares, and thus knowing how many such small structures are likely to appear is critical. One simple example of a small substructure in a Latin square is an \emph{intercalate}, where two rows and two columns intersect on a 2 by 2 subsquare which is itself a Latin square (see Figure~\ref{fig:LS} for an example). Using switching arguments, Kwan and Sudakov~\cite{kwan2018intercalates} showed that a random Latin square $L_n$ has at least $(1-o(1))\frac{n^2}{4}$ intercalates with high probability. This lower bound matched a conjecture of McKay and Wanless~\cite{mckay1999most}, whose proof was recently completed by Kwan, Sah and Sawhney~\cite{kwan2022large} showing that a random Latin square $L_n$ has at most $(1+o(1))\frac{n^2}{4}$ intercalates with high probability.
 If Conjecture~\ref{conj:random} is to be proved by an absorption approach building on those in \cite{gould2023hamilton,kwan2020almost}, much of a random Latin square will need to be decomposed into small structures which are used to build absorbers. In recent years, related decompositions have been found by applying hypergraph matching results (see Ehard, Glock and Joos~\cite{ehard2020pseudorandom} for such a result) to some almost-regular auxiliary hypergraphs in an increasingly sophisticated way.
 For such an approach, upper bounds for the number of specific small substructures are needed as well as lower bounds, and thus the work of Kwan, Sah and Sawhney~\cite{kwan2022large} (along with  further developments by the same authors and Simkin~\cite{kwan2022substructures}) is important here, where the upper bound is proved using the switching methods pioneered by McKay and Wanless~\cite{mckay1999most} and the deletion method of R\"odl and Ruci\'nski~\cite{janson2002infamous,rodl1995threshold}.

%%%%%%%%%%%%%%%%%%%%%%%%%%%%%%%%%%%%%%%%%%%%%%%%%%%%%%%%%%%%%%%%%%%%%%%%%%%%%%%%%%%%%%%%%%%%%%%%%%%%%%%%%%%%%%%%%%%%%%%%%%
%%%%%%%%%%%%%%%%%%%%%%%%%%%%%%%%%%%%%%%%%%%%%%%%%%%%%%%%%%%%%%%%%%%%%%%%%%%%%%%%%%%%%%%%%%%%%%%%%%%%%%%%%%%%%%%%%%%%%%%%%%
%%%%%%%%%%%%%%%%%%%%%%%%%%%%%%%%%%%%%%%%%%%%%%%%%%%%%%%%%%%%%%%%%%%%%%%%%%%%%%%%%%%%%%%%%%%%%%%%%%%%%%%%%%%%%%%%%%%%%%%%%%
%%%%%%%%%%%%%%%%%%%%%%%%%%%%%%%%%%%%%%%%%%%%%%%%%%%%%%%%%%%%%%%%%%%%%%%%%%%%%%%%%%%%%%%%%%%%%%%%%%%%%%%%%%%%%%%%%%%%%%%%%%

\section{Latin squares without full transversals and related conjectures}\label{sec:examples}
\begin{figure}
\begin{center}
\begin{tikzpicture}[define rgb/.code={\definecolor{mycolor}{rgb}{#1}}, rgb color/.style={define rgb={#1},mycolor},scale=1.2]

\def\wi{0.4cm}

\def\n{5}

\def\Latentry{{{A,K,Q,J}, {Q,J,A,K},{J,Q,K,A},{K,A,J,Q}}};

\foreach \x in {0,1,...,3}
\foreach \y in {0,1,...,3}
{
\coordinate (A\x\y) at ($\n*(0,\wi)+\x*(\wi,0)-\y*(0,\wi)$);
\draw (A\x\y) [thick] rectangle ($(A\x\y)+(\wi,\wi)$);
}

\foreach \x/\y in {1/0,3/1,0/2,2/3}
{
%\draw (A\x\y) [thick,fill=black!40] rectangle ($(A\x\y)+(\wi,\wi)$);
}

\draw ($(A10)+0.5*(\wi,\wi)$) [thick,densely dotted] circle [radius = 0.45*\wi];% ($(A10)+(\wi,\wi)$);
\draw ($(A31)+0.5*(\wi,\wi)$) [thick,black!25] circle [radius = 0.45*\wi];% ($(A31)+(\wi,\wi)$);
\draw ($(A02)+0.5*(\wi,\wi)$) [thick] circle [radius = 0.45*\wi];% ($(A02)+(\wi,\wi)$);
\draw ($(A23)+0.5*(\wi,\wi)$) [thick,densely dashed] circle [radius = 0.45*\wi];% ($(A23)+(\wi,\wi)$);

\foreach \x in {0,1,...,3}
\foreach \y in {0,1,...,3}
{
%\coordinate (A\x\y) at ($\n*(0,\wi)+\x*(\wi,0)-\y*(0,\wi)$);
\draw (A\x\y) [thick] rectangle ($(A\x\y)+(\wi,\wi)$);
}

\foreach \x/\y in {1/0,2/1,3/2,0/3}
{
\draw ($(A\x\y)+0.5*(\wi,\wi)$) node {1};
}
\foreach \x/\y in {0/0,3/1,2/2,1/3}
{
\draw ($(A\x\y)+0.5*(\wi,\wi)$) node {2};
}
\foreach \x/\y in {3/0,0/1,1/2,2/3}
{
\draw ($(A\x\y)+0.5*(\wi,\wi)$) node {3};
}
\foreach \x/\y in {2/0,1/1,0/2,3/3}
{
\draw ($(A\x\y)+0.5*(\wi,\wi)$) node {4};
}
\foreach \x/\y in {3/0,0/1,4/2,1/3,2/4,5/5}
{
%\draw ($(A\x\y)+0.5*(\wi,\wi)$) node {5};
}
\foreach \x/\y in {2/0,5/1,0/2,3/3,1/4,4/5}
{
%\draw ($(A\x\y)+0.5*(\wi,\wi)$) node {6};
}

\end{tikzpicture}\hspace{0.5cm}\begin{tikzpicture}\draw [white] (0,0) -- (0,-0.9); \draw (0,0) node {$\implies$};\end{tikzpicture}\hspace{0.5cm}\begin{tikzpicture}[define rgb/.code={\definecolor{mycolor}{rgb}{#1}}, rgb color/.style={define rgb={#1},mycolor},scale=1.2]

\def\wi{0.4cm}

\def\n{5}
\def\vx{0.05cm}

\def\Latentry{{{A,K,Q,J}, {Q,J,A,K},{J,Q,K,A},{K,A,J,Q}}};

\foreach \x in {0,1}
\foreach \y in {0,1,2,3}
{
\coordinate (A\x\y) at ($\n*(0,\wi)+2*\x*(\wi,0)-\y*(0,\wi)+0.5*(\wi,\wi)$);
%\draw [fill] ($(A\x\y)+0.5*(\wi,\wi)$) circle [radius=\vx];
%\draw [thick] rectangle ($(A\x\y)+(\wi,\wi)$);
}
%\foreach\x in {0,1}

\foreach \y/\x in {0/1,1/2,2/3,3/4}
{
\draw ($(A0\y)-(0.25,0)$) node {\footnotesize \x};
\draw ($(A1\y)+(0.25,0)$) node {\footnotesize \x};
}

\foreach \x/\y in {1/0,3/1,0/2,2/3}
{
%\draw (A\x\y) [thick,fill=black!40] rectangle ($(A\x\y)+(\wi,\wi)$);
}

\foreach \x/\y in {1/0,2/1,3/2,0/3}
{
\draw [thick,densely dashed] (A0\x) -- (A1\y);
}
\foreach \x/\y in {0/0,3/1,2/2,1/3}%correct
{
\draw [thick,black!25] (A0\x) -- (A1\y);
}
\foreach \x/\y in {3/0,0/1,1/2,2/3}
{
\draw [thick,black,densely dotted] (A0\x) -- (A1\y);
}
\foreach \x/\y in {2/0,1/1,0/2,3/3}
{
\draw [thick] (A0\x) -- (A1\y);
}

\foreach \x in {0,1}
\foreach \y in {0,1,2,3}
{
%\coordinate (A\x\y) at ($\n*(0,\wi)+2*\x*(\wi,0)-\y*(0,\wi)+0.5*(\wi,\wi)$);
\draw [fill] (A\x\y) circle [radius=\vx];
%\draw [thick] rectangle ($(A\x\y)+(\wi,\wi)$);
}

\draw ($(A00)+(-0.25,0.3)$) node {\footnotesize row};
\draw ($(A10)+(0.25,0.35)$) node {\footnotesize column};

\end{tikzpicture}\hspace{1.2cm}\begin{tikzpicture}\draw [thick,dotted] (0,1.1) -- (0,-1.1);\end{tikzpicture}\hspace{1.2cm}\begin{tikzpicture}[define rgb/.code={\definecolor{mycolor}{rgb}{#1}}, rgb color/.style={define rgb={#1},mycolor},scale=1.2]

\def\wi{0.4cm}

\def\n{5}
\def\vx{0.05cm}

\def\Latentry{{{A,K,Q,J}, {Q,J,A,K},{J,Q,K,A},{K,A,J,Q}}};

\foreach \x in {0,1}
\foreach \y in {0,1,2,3}
{
\coordinate (A\x\y) at ($\n*(0,\wi)+2*\x*(\wi,0)-\y*(0,\wi)+0.5*(\wi,\wi)$);
%\draw [fill] ($(A\x\y)+0.5*(\wi,\wi)$) circle [radius=\vx];
%\draw [thick] rectangle ($(A\x\y)+(\wi,\wi)$);
}
%\foreach\x in {0,1}

\foreach \y/\x in {0/1,1/2,2/3,3/4}
{
\draw ($(A0\y)-(0.2,0)$) node {\footnotesize \x};
\draw ($(A1\y)+(0.2,0)$) node {\footnotesize \x};
}

\foreach \x/\y in {1/0,3/1,0/2,2/3}
{
%\draw (A\x\y) [thick,fill=black!40] rectangle ($(A\x\y)+(\wi,\wi)$);
}

\foreach \x/\y in {0/1}
{
\draw [thick,densely dotted] (A0\x) -- (A1\y);
\draw [black!50] ($0.5*(A0\x)+0.5*(A1\y)+(-0.1,0.195)$) node {\footnotesize 1};
}
\foreach \x/\y in {1/3}
{
\draw [thick,black!25] (A0\x) -- (A1\y);
\draw [black!50] ($0.5*(A0\x)+0.5*(A1\y)+(-0.1,-0.05)$) node {\footnotesize 2};
}
\foreach \x/\y in {3/2}
{
\draw [thick,black,densely dashed] (A0\x) -- (A1\y);
\draw [black!50] ($0.5*(A0\x)+0.5*(A1\y)-(0,0.15)$) node {\footnotesize 3};
}
\foreach \x/\y in {2/0}
{
\draw [thick] (A0\x) -- (A1\y);
\draw [black!50] ($0.5*(A0\x)+0.5*(A1\y)+(0.05,-0.075)$) node {\footnotesize 4};
}

\foreach \x in {0,1}
\foreach \y in {0,1,2,3}
{
%\coordinate (A\x\y) at ($\n*(0,\wi)+2*\x*(\wi,0)-\y*(0,\wi)+0.5*(\wi,\wi)$);
\draw [fill] (A\x\y) circle [radius=\vx];
%\draw [thick] rectangle ($(A\x\y)+(\wi,\wi)$);
}
\draw ($(A00)+(-0.25,0.3)$) node {\footnotesize \textcolor{black}{row}};
\draw ($(A10)+(0.25,0.35)$) node {\footnotesize \textcolor{black}{column}};
\end{tikzpicture}
\end{center}
\caption{A Latin square of order 4 converted into an edge-coloured complete bipartite graph. A transversal in the Latin square is shown by highlighting each cell with the style used to represent the edges with colour corresponding to the cell entry. The corresponding perfect rainbow matching is then shown on the right.}\label{fig:translate}
\end{figure}
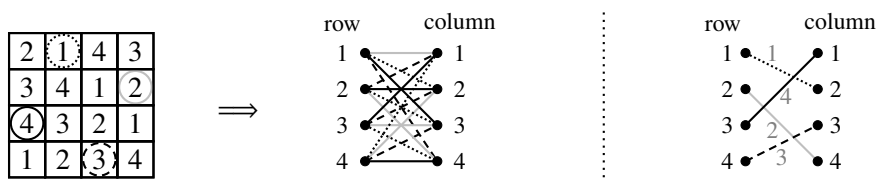

It is convenient, and now commonplace, to translate transversals in Latin squares into a rainbow matching problem in optimally coloured complete bipartite balanced graphs, as follows (see also Figure~\ref{fig:translate}). Given a Latin square $L$, create a vertex corresponding to each row and each column, and, for each row/column pair put an edge between their corresponding vertices coloured by the symbol in the cell indexed by that row and column. This gives a coloured bipartite graph $H(L)$, which is a copy of the complete bipartite $n$ by $n$ graph $K_{n,n}$ where the edges are coloured using the set of symbols of the Latin square. That each symbol appears only once in each row and in each column of the Latin square implies that every colour appears at each vertex at most once, and therefore the colouring is \emph{proper}. As the colouring uses the minimal number of colours for a proper colouring of $K_{n,n}$, we say it is \emph{optimally coloured}. Given any optimal colouring of $K_{n,n}$, the reverse construction easily creates a corresponding Latin square of order $n$, and therefore the Latin squares of order $n$ correspond exactly to the optimal colourings of the complete bipartite graph $K_{n,n}$.

Any transversal in a Latin square $L$ corresponds to a subgraph in the corresponding graph $H(L)$ in which no edges share any vertices or any colours (see Figure~\ref{fig:translate}). Subgraphs in edge-coloured graphs in which each colour appears at most once are known as \emph{rainbow}, where rainbow subgraphs have in recent years seen a hive of activity (see, for example, the many results covered by the survey by Pokrovskiy~\cite{Alexeysurvey}). A partial transversal then in $L$ corresponds to a rainbow matching in $H(L)$, while a full transversal in $L$ corresponds to a perfect rainbow matching in $H(L)$.

This translation allows us to use graph theory notation and techniques in the study of Latin squares. For example, for any abelian group $G$ of order $n$, let $L(G)$ be the Latin square corresponding to its addition table and consider the graph, $H(G)$ say, corresponding to $L(G)$. If $H(G)$ contains a perfect rainbow matching, then, letting $c(e)$ be the colour of the edge $e$, we have
\begin{equation}\label{eqn:sum}
\sum_{v\in G}v=\sum_{e\in M}c(e)=\sum_{xy\in M}(x+y)=2\sum_{v\in G}v,
\end{equation}
so that we must have that $\sum_{v\in G}v=0$, where 0 is the identity of $G$. Where this does not hold, then, $H(G)$ can have no perfect rainbow matching. Where $G=\mathbb{Z}_{2m}$, for any integer $m\geq 1$, this corresponds to the example given in Section~\ref{sec:intro}, but already this gives us more examples of Latin squares without a full transversal, for example by taking $G=\mathbb{Z}_{6}\times \mathbb{Z}_3$.

In fact, as, for simplicity, we have set $G$ to be abelian, a characterisation of when the addition table for $G$ has a full transversal has long been known due to Paige~\cite{paige1947note} in 1947 (using the language of complete mappings). That is, there is not a full transversal exactly when $G$ has exactly one element of order 2, or, equivalently, when its only Sylow 2-subgroup is a cyclic group. In this case, this condition can be easily seen to be equivalent to $\sum_{v\in G}v\neq 0$.

In 1950, Paige~\cite{paige1951complete} then considered whether $L(G)$ has a full transversal or not for finite non-abelian groups $G$, where $L(G)$ is the Latin square corresponding to the multiplication table of $G$. Paige proved that if $L(G)$ has a full transversal then there must be an ordering of the elements of $G$ such that their product in that ordering is equal to the identity. Noting that this was shown to be a sufficient condition for finite abelian groups in~\cite{paige1947note}, Paige conjectured that, likewise, this is a sufficient condition for a full transversal to occur in the multiplication table of a finite non-abelian group. As Paige noted, this is easy to show for finite non-abelian groups of odd order as the leading diagonal of the corresponding Latin squares forms a full transversal.

In 1955, Hall and Paige~\cite{hallpaige} considered an equivalent condition, that $\sum_{x\in G}x=0$ in the abelianisation $G^{\mathrm{ab}}=G/G'$ of $G$, where $G'$ is the commutator subgroup of $G$. This condition is now known as the Hall-Paige condition. Hall and Paige~\cite{hallpaige} conjectured that this is sufficient for a full transversal to occur in the multiplication table for any finite group, and also showed that this condition is equivalent to the condition that all the Sylow 2-subgroups of $G$ are trivial or non-cyclic.
The Hall-Paige conjecture was finally confirmed in 2009 through a combination of work by Wilcox, Evans and Bray~\cite{wilcox,evans}, as described in Section~\ref{sec:full}, and recorded below.

\begin{conj}[The Hall-Paige conjecture, proved by Wilcox, Evans and Bray]\label{conj:HP} For any finite group $G$, $L(G)$ has a full transversal if and only if $G$ satisfies the Hall-Paige condition.% $\sum_{x\in G}x=0$ in $G^{\mathrm{ab}}$.
\end{conj}

In the context of Section~\ref{sec:decomp}, we note that, for any finite group $G$, if $L(G)$ has a full transversal then it can be very easily seen to have a decomposition into full transversals. Indeed, if $f:G\to G$ is a bijection such that the cells indexed by $(g,f(g))$, $g\in G$, form a full transversal then, for each $h\in G$, the cells indexed by $(g,f(g)\circ h)$, $g\in G$, are also a full transversal, and thus gives a decomposition into full transversals.

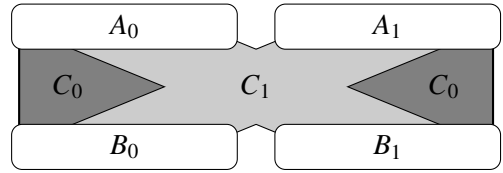
\begin{figure}
\begin{center}
\begin{tikzpicture}[scale=1]
\def\vxrad{0.07cm}
\def\horunit{0.8}
\def\verunit{1.75}
\def\edgelength{0.4}
\def\betweenrows{0.5}

%\draw [white] ($-1*(\horunit,0)$) -- ($41*(\horunit,0)$);

\def\gapratio{1}

\foreach \lab/\yy in {A/1,B/-1}
{
\foreach \num/\xx in {1/-1,2/1}
{
\coordinate (\lab\num) at ($(\xx*\verunit,\yy*\horunit)$);
}
}

\coordinate (L11) at ($(A1)-(1,0)-(0.2,0)$);
\coordinate (L12) at ($(B1)-(1,0)-(0.2,0)$);
\coordinate (L21) at ($(A2)-(1,0)$);
\coordinate (L22) at ($(B2)-(1,0)$);
\coordinate (R11) at ($(A1)+(1,0)$);
\coordinate (R12) at ($(B1)+(1,0)$);
\coordinate (R21) at ($(A2)+(1,0)+(0.2,0)$);
\coordinate (R22) at ($(B2)+(1,0)+(0.2,0)$);

\foreach \cooo in {L11,L12}
\foreach \w in {w}
{
\coordinate (\cooo\w) at ($(\cooo)-(0.2,0)$);
}
\foreach \cooo in {L21,L22}
\foreach \w in {w}
{
\coordinate (\cooo\w) at ($(\cooo)-(0.45,0)$);
}
\foreach \cooo in {R21,R22}
\foreach \w in {w}
{
\coordinate (\cooo\w) at ($(\cooo)+(0.2,0)$);
}
\foreach \cooo in {R11,R12}
\foreach \w in {w}
{
\coordinate (\cooo\w) at ($(\cooo)+(0.45,0)$);
}

\fill [fill=black!50] ($(L11w)$) -- ($(L12w)$) -- (R12w) -- (R11w) -- (L11w);
\fill [fill=black!50] ($(L21w)$) -- ($(L22w)$) -- (R22w) -- (R21w) -- (L21w);
\draw [thick] ($(L11w)$) -- ($(L12w)$);
\draw [thick] ($(L21w)$) -- ($(L22w)$);
\draw [thick] ($(R21w)$) -- ($(R22w)$);
\draw [thick] ($(R11w)$) -- ($(R12w)$);

\draw [black,thick] (L11w) -- (L22) -- (R22w) -- (R11) -- (L11w);
\draw [black,thick]  (L21) -- (L12w) -- (R12) -- (R21w) -- (L21);
\fill [black,fill=black!20] (L11w) -- (L22) -- (R22w) -- (R11) -- (L11w);
\fill [black,fill=black!20] (L21) -- (L12w) -- (R12) -- (R21w) -- (L21w);
%\draw [red] (L11w) -- (L22) -- (R22w) -- (R11) -- (L11w);
%\draw [red]  (L21) -- (L12w) -- (R12) -- (R21w) -- (L21);

\foreach \lab/\num/\labb in {A/1/A_0,A/2/A_1,B/1/B_0,B/2/B_1}
{
%\draw [fill=white] (\lab\num) circle [x radius=1.5cm,y radius=0.3cm];
\draw [rounded corners,fill=white] ($(\lab\num)-(1.5,0.3)$) rectangle ($(\lab\num)+(1.5,0.3)$) {};
\draw (\lab\num) node {$\labb$};
}

\draw  ($0.5*(A1)+0.5*(B1)-(0.75,0)$) node {$C_0$};
\draw  ($0.5*(A2)+0.5*(B2)+(0.75,0)$) node {$C_0$};
\draw ($0.25*(A1)+0.25*(B1)+0.25*(A2)+0.25*(B2)$) node {$C_1$};

\end{tikzpicture}
\end{center}
\caption{A blow-up construction based on the addition table for $\mathbb{Z}_2$. For some $m$, each labelled set has size $m$. The edges between $A_0$ and $B_0$, and between $A_1$ and $B_1$, are all present and properly coloured using the colours in $C_0$, while the edges between $A_0$ and $B_1$, and between $A_1$ and $B_0$, are all present and properly coloured using colours in $C_1$. When $m$ is odd, there is no perfect rainbow matching.}\label{fig:moreeg}
\end{figure}

Given a finite group $G$ which does not satisfy the Hall-Paige condition, we can generate further Latin squares by using what, in modern parlance, is a `blow-up' construction (see Figure~\ref{fig:moreeg} for a depiction with $G=\mathbb{Z}_2$). For convenience we will assume that $G$ is an abelian group, of order $n$, and construct optimal colourings of complete bipartite graphs. Let $k\geq 1$ be an integer. Disjointly, for each $v\in G$, take sets $A_v$ and $B_v$, each of $k$ vertices, and a set $C_v$ of $k$ colours. Let $A=\cup_{v\in G}A_v$ and $B=\cup_{v\in G}B_v$ and, for each pair $v$ and $w$, put all possible edges between $A_v$ and $B_w$ and properly colour them using the colours in $C_{v+w}$.
This gives a graph $H$ which is an optimally coloured copy of $K_{kn,kn}$ using the colours in $C=\cup_{v\in G}C_v$. If $H$ contains a perfect rainbow matching $M$, then, for each edge $e\in M$, let $v_e,w_e\in G$ be such that $e$ lies between $A_{v_e}$ and $B_{w_e}$, noting that $c(e)\in C_{v_e+w_e}$. Then, the corresponding version of \eqref{eqn:sum} is
\begin{equation*}
k\cdot \sum_{u\in G}u=\sum_{u\in G}|C_u|\cdot u=\sum_{u\in G}\sum_{e\in M:c(e)\in C_u}u=\sum_{e\in M}(v_e+w_e)
=2k\cdot \sum_{u\in G}u.
%|\cdot v=\sum_{}\sum_{xy\in M:}c(xy)=\sum_{xy\in M}(x+y)=2\sum_{v\in G}v,k\cdot \sum_{v\in G}v=\sum_{v\in G}\sum_{c\in C_v}v=\sum_{c\in C}\sum_{xy\in M:c(xy)\in C_v}
%=k\cdot \sum_{v\in G}v+k\cdot \sum_{w\in G}w=
\end{equation*}
Therefore, if $k\cdot \sum_{u\in G}u\neq 0$, then $H$ can have no perfect rainbow matching, and thus the corresponding Latin square has no full transversal.

In the case $G=\mathbb{Z}_{2m}$ for any $m\geq 1$ and odd $k$, this example was given by Maillet~\cite{maillet1894carres} in 1894 (and later rediscovered by Parker~\cite{parker1971pathological}).
As observed by Cavenagh and Wanless~\cite{cavenagh2016latin}, this case already provides many Latin squares with no full transversals. For example, when $n=2k$ and $k$ is odd then the number of Latin squares with no full transversals is at least $\big(\big(\frac{1}{2}e^{-2}-o(1)\big)n\big)^{n^2}$. Of course, as shown by the result of Kwan~\cite{kwan2020almost} discussed in Section~\ref{sec:decomp}, this is a vanishingly small proportion of Latin squares of order $n$, as the number of Latin squares of order $n$ is $|\mathcal{L}(n)|=((e^{-2}+o(1))n)^{n^2}$ (see~\cite[Chapter 17]{van2001course}).

All the examples of Latin squares without a full transversal that we have covered here have even order. In 1967, Ryser~\cite{Ryser} (see also~\cite{Rysertranslation}) conjectured that there are no Latin squares of odd order without a full transversal. Brualdi~\cite{Brualdi} later conjectured that every Latin square of order $n$ has a partial transversal with $n-1$ cells, while Stein~\cite{stein} made some related conjectures in 1975 which are stronger than this (see Section~\ref{sec:related}). For at least the last decade, the following well-known combined conjecture has been known as the Ryser-Brualdi-Stein conjecture.

\begin{conj}[The Ryser-Brualdi-Stein conjecture]\label{conj:RBS}
Every Latin square of order $n$ has a partial transversal with at least $n-1$ cells, and a full transversal if $n$ is odd.
\end{conj}

In Section~\ref{sec:full} we will discuss the proof of the Hall-Paige conjecture and other cases where Latin squares are known to have full transversals, before discussing progress towards the Ryser-Brualdi-Stein conjecture in Section~\ref{sec:partial}. Before this, however,  let us briefly note that we have seen no Latin squares of odd order $n$ with any meaningful restriction on their transversal properties. However, there are Latin squares of order $n$ which have no decomposition into full transversals for each odd $n\geq 5$. This was shown for every $n\geq 5$ with $n\equiv 1\;\mathrm{mod}\; 4$ by Mann~\cite{mann1944orthogonal}, and for every $n\geq 7$ with $n\equiv 3\;\mathrm{mod}\; 4$ independently by Evans~\cite{evans2006latin} and Wanless and Webb~\cite{wanless2006existence}.% for every odd $n$ except $n=1,3$ there is some Latin square of order $n$ which has no decomposition into full transversals.

%In the rest of this section we will discuss further examples of Latin squares without full transversals. While doing so we will state several conjectures (most notably Conjecture~\ref{conj:RBS}), but postpone any discussion of the progress towards them to Section~\ref{sec:full} and Section~\ref{sec:partial}.

%%%%%%%%%%%%%%%%%%%%%%%%%%%%%%%%%%%%%%%%%%%%%%%%%%%%%%%%%%%%%%%%%%%%%%%%%%%%%%%%%%%%%%%%%%%%%%%%%%%%%%%%%%%%%%%%%%%%%%%%%%
%%%%%%%%%%%%%%%%%%%%%%%%%%%%%%%%%%%%%%%%%%%%%%%%%%%%%%%%%%%%%%%%%%%%%%%%%%%%%%%%%%%%%%%%%%%%%%%%%%%%%%%%%%%%%%%%%%%%%%%%%%
%%%%%%%%%%%%%%%%%%%%%%%%%%%%%%%%%%%%%%%%%%%%%%%%%%%%%%%%%%%%%%%%%%%%%%%%%%%%%%%%%%%%%%%%%%%%%%%%%%%%%%%%%%%%%%%%%%%%%%%%%%
%%%%%%%%%%%%%%%%%%%%%%%%%%%%%%%%%%%%%%%%%%%%%%%%%%%%%%%%%%%%%%%%%%%%%%%%%%%%%%%%%%%%%%%%%%%%%%%%%%%%%%%%%%%%%%%%%%%%%%%%%%

\section{Full transversals in Latin squares}\label{sec:full}
From its statement in 1955 until its proof in 2009, the Hall-Paige conjecture (Conjecture~\ref{conj:HP}) was steadily shown to hold for various groups and saw much related work. For a detailed exposition of this, as well as a complete unified proof of the conjecture, we refer the reader to Part II of the book by Evans~\cite{evans2018orthogonal}. In their original paper, Hall and Paige~\cite{hallpaige} showed that, if $G$ has a normal subgroup $N$ for which both $L(N)$ and $L(G/N)$ has a full transversal, then $L(G)$ has a full transversal, and used this to prove Conjecture~\ref{conj:HP} for solvable groups.
In a major breakthrough in 2009, also via an inductive argument, Wilcox~\cite{wilcox} reduced Conjecture~\ref{conj:HP} to the case of simple groups. Thus, the classification of finite simple groups exactly characterised the remaining task. Hall and Paige~\cite{hallpaige} had already confirmed the conjecture for alternating groups. Wilcox~\cite{wilcox} gave a unified proof for groups of Lie type, while Evans~\cite{evans} combined Wilcox's methods with computer algebra to give a proof for the Tits group and all but one of the 26 sporadic groups. Bray then checked the final remaining case, which was the fourth Janko group (see~\cite{wilcox}, and \cite{bray2020hall} for the eventual publication), which completed the proof of the Hall-Paige conjecture.

More recently, Eberhard, Manners and Mrazovi{\'c} \cite{greenalites} gave a completely different proof of the Hall-Paige conjecture for large groups, using tools from analytic number theory. Remarkably, they were able to show that any group $G$ satisfying the Hall-Paige condition has
$\big(e^{-1/2}+o(1)\big)|G^{\mathrm{ab}}|\frac{(n!)^2}{n^n}$ full transversals, where, again, $G^{\mathrm{ab}}$ is the abelianisation of $G$, thus giving a precise asymptotic. (Indeed, they even determine the order of the $o(1)$ term, see~\cite[Theorem 1.4]{greenalites}!)
The quantitative bounds in Eberhard, Manners and Mrazovi{\'c}'s work are strong enough to already rule out many of the cases considered in the original proof of the Hall-Paige conjecture by Wilcox, Bray and Evans. This gives a proof of Conjecture~\ref{conj:HP} which avoids extensive case-checking for the sporadic groups or the Tits group, from these only requiring verification of the conjecture for the first two Mathieu groups (see~\cite{greenalites}).

Even more recently, Eberhard, Manners and Mrazovi{\'c}~\cite{eberhard2023transversals} again used tools from analytic number theory (including a loose variant of the circle method) to find full transversals in Latin squares satisfying a quasirandomness condition. Here the quasirandomness condition is defined in terms of the spectral gap of an operator associated with the Latin square (see~\cite{eberhard2023transversals}).
Where the Latin square coincides with the multiplication table of a group, their condition coincides with the definition of quasirandomness for subgroups by
Gowers~\cite{gowers2008quasirandom}, and recovers the main result of their work in~\cite{greenalites} for sufficiently quasirandom groups. Using recent results of Kwan, Sah, Sawhney and Simkin~\cite{kwan2022substructures}, Eberhard, Manners and Mrazovi\'c showed that a random Latin square of order $n$ satisfies their quasirandomness condition with high probability, allowing them to improve Kwan's result (as discussed in Section~\ref{sec:decomp}) to show that a random Latin square of order $n$ has $\big(e^{-1/2}+o(1)\big)\frac{(n!)^2}{n^n}$ full transversals with high probability.

A further independent proof of the Hall-Paige conjecture has been also recently been given, this time with  combinatorial methods, by M\"uyesser and Pokrovskiy~\cite{muyesser2022random} as part of a generalised theorem with several applications resolving old problems in combinatorial group theory in the case of large groups (see~\cite{muyesser2022random}). To get an idea of this general theorem, first let $G$ be a large group of order $n$ and let $H(G)$ be the bipartite graph corresponding to the multiplication table of $G$ (as described in Section~\ref{sec:examples}) with vertex sets $A$ and $B$ and colour set $C$. Suppose we have sets $A'\subset A$, $B'\subset B$ and $C'\subset C$ of equal size.
When can we expect a perfect rainbow matching between $A'$ and $B'$ in $H(G)$ using only the colours in $C'$? Similarly to the arguments we have covered, it is easy to see that a necessary condition must be that $\sum_{a\in A'}a+\sum_{b\in B'}b=\sum_{c\in C'}c$ in the abelianisation $G^{\mathrm{ab}}$ of $G$. M\"uyesser and Pokrovskiy~\cite{muyesser2022random} show that this is a sufficient condition for sets $A'$, $B'$ and $C'$ which are sufficiently close to random subsets where each element is independently included with the same probability $p$, for any $p\geq n^{-1/10^{100}}$. Taking $p=1$ in this result, then, recovers a proof of the Hall-Paige conjecture for large groups.

The main tool used by M\"uyesser and Pokrovskiy in~\cite{muyesser2022random} is the absorption method, as described loosely in Section~\ref{sec:decomp}. As we will discuss in detail an absorption approach to finding large rainbow matchings in the next section, we will outline how this can be used now. %2
(This discussion is representative of the methods in~\cite{muyesser2022random}, but slightly different to align with our subsequent discussion.)
A natural approach for absorption here would be, for some $\ell_0,\ell_1\leq n$, to look for a balanced set $V^{\mathrm{abs}}$ of $2\ell_0$ vertices and a set $C^{\mathrm{abs}}$ of $\ell_0+\ell_1$ colours such that, for any balanced set $W\subset V(H)\setminus V^{\mathrm{abs}}$ of $2\ell_1$ vertices, there is a rainbow matching in $H:=H(G)$ with vertex set $V^{\mathrm{abs}}\cup W$ and colour set $C^{\mathrm{abs}}$. (A set is \emph{balanced} if it has equally many vertices in each side of the bipartition.)
In the terminology of absorption, we would say that $(V^{\mathrm{abs}},C^{\mathrm{abs}})$ can absorb any balanced vertex set with $2\ell_1$ vertices. If this could be done, then the strategy would look for an initial rainbow matching with vertex set in $V(G)\setminus V^{\mathrm{abs}}$ which uses exactly the colours not in $C^{\mathrm{abs}}$ (made easier than the original problem as there are $2\ell_1$ more vertices than have to be used in this matching) before applying the absorption property to the unused vertices in $V(H)\setminus V^{\mathrm{abs}}$ to find a second matching which, added to the initial matching, completes a perfect rainbow matching in $H$.

However, there is an inherent complication here. If there is a rainbow matching in $H$ with vertex set $V^{\mathrm{abs}}\cup W$ and colour set $C^{\mathrm{abs}}$, then, assuming from now for simplicity that $G$ is abelian, we have
\begin{equation}\label{eqn:recall}
\sum_{v\in V^{\mathrm{abs}}\cup W}v=\sum_{c\in C^{\mathrm{abs}}}c,
\end{equation}
and therefore the best we can hope for is that $(V^{\mathrm{abs}},C^{\mathrm{abs}})$ can absorb any balanced set of $2\ell_1$ vertices for which $\sum_{v\in W}v=\sum_{c\in C^{\mathrm{abs}}}-\sum_{v\in V^{\mathrm{abs}}}v$.
We may as well say this sum is equal to 0, and hope to create an absorption structure $(V^{\mathrm{abs}},C^{\mathrm{abs}})$ which can absorb any appropriate set of vertices with zero sum.

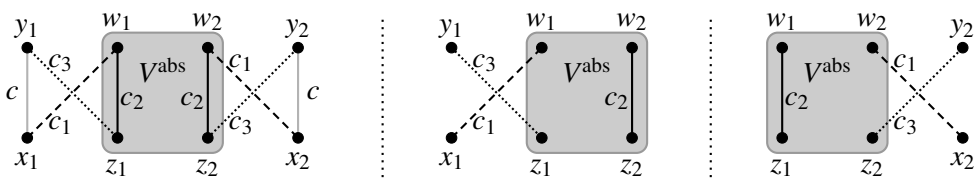
\begin{figure}[b]
\begin{center}
\begin{tikzpicture}[scale=1]
\def\vxrad{0.07cm}
\def\horunit{1.2}
\def\edgelength{0.4}
\def\betweenrows{0.5}

%\draw [white] ($-1*(\horunit,0)$) -- ($41*(\horunit,0)$);

\def\gapratio{1}

\foreach \num/\parity/\parityy in {1/-1/-1,2/-1/1}
{
\coordinate (Y\num) at ($(0,0)+\parityy*\gapratio*(\horunit,0)+(0.5*\parity*\horunit,0.5*\horunit)$);
\coordinate (W\num) at ($(0,0)+\parityy*\gapratio*(\horunit,0)+(-0.5*\parity*\horunit,0.5*\horunit)$);
\coordinate (X\num) at ($(0,0)+\parityy*\gapratio*(\horunit,0)+(0.5*\parity*\horunit,-0.5*\horunit)$);
\coordinate (Z\num) at ($(0,0)+\parityy*\gapratio*(\horunit,0)+(-0.5*\parity*\horunit,-0.5*\horunit)$);
}

\def\drop{0.2};
\draw [rounded corners,thick,black!40,fill=black!20] ($0.5*(Z1)+0.5*(X2)-(0,\drop)$)
-- ($(X2)+(\drop,-\drop)$) -- ($(Y2)+(\drop,\drop)$) -- ($(W1)+(-\drop,\drop)$) -- ($(Z1)+(-\drop,-\drop)$) -- ($0.5*(Z1)+0.5*(X2)-(0,\drop)$);
\draw ($0.25*(W1)+0.25*(Z1)+0.25*(Y2)+0.25*(X2)+(0,0.3)$) node {$V^{\mathrm{abs}}$};

%%%EDGES
\foreach \x/\y/\col/\n in {Z/W/teal/1,X/Y/teal/2}
{
\draw [thick] (\x\n) -- (\y\n);
}

\foreach \x/\y/\col in {X/W/red}
\foreach \n in {1}
{
\draw [thick,densely dashed] (\x\n) -- (\y\n);
}

\foreach \x/\y/\col in {Y/Z/red}
\foreach \n in {2}
{
\draw [thick,densely dashed] (\x\n) -- (\y\n);
}

\foreach \x/\y/\col in {Y/Z/blue}
\foreach \n in {1}
{
\draw [thick,densely dotted] (\x\n) -- (\y\n);
}
\foreach \x/\y/\col in {X/W/blue}
\foreach \n in {2}
{
\draw [thick,densely dotted] (\x\n) -- (\y\n);
}

\foreach \x/\y/\col in {X/Y/magenta}
\foreach \n in {1}
{
\draw [thick,black!30] (\x\n) -- (\y\n);
}
\foreach \x/\y/\col in {W/Z/magenta}
\foreach \n in {2}
{
\draw [thick,black!30] (\x\n) -- (\y\n);
}

%%vertex labels
\foreach \x/\lab/\offs/\sss in {X/x/-1.2/5,Y/y/1/6}%,Z/z/-1.2/7,W/w/1/4}
\foreach \num in {1}
{
\draw  ($(\x\num)+\offs*(0,0.25)$) node {$\lab_\num$};
}
%%%%%%%%%%%%%%%NEW LABELS%%%%%%%%%%%%%%%%%%%
\foreach \x/\lab/\offs/\sss in {Z/z/-1.6/5,W/w/1.4/4}
\foreach \num in {1}
{
\draw  ($(\x\num)+\offs*(0,0.25)$) node {$\lab_\num$};
}
\foreach \x/\lab/\offs/\sss in {X/z/-1.6/5,Y/w/1.4/6}
\foreach \num in {2}
{
\draw  ($(\x\num)+\offs*(0,0.25)$) node {$\lab_\num$};
}
%%%%%%%%%%%%%%%NEW LABELS%%%%%%%%%%%%%%%%%%%

\foreach \x/\lab/\offs/\sss in {Z/x/-1.2/7,W/y/1/4}%X/z/-1.2/5,Y/w/1/6,
\foreach \num in {2}
{
\draw  ($(\x\num)+\offs*(0,0.25)$) node {$\lab_\num$};
}

%%%colourlabels

\draw ($0.5*(X2)+0.5*(W2)+(0.025,-0.25)+0.15*(-\horunit,-\horunit)$) node {$c_3$};
\draw ($0.5*(X1)+0.5*(W1)+(0.025,-0.25)+0.15*(-\horunit,-\horunit)$) node {$c_1$};

\draw ($0.5*(X2)+0.5*(Y2)-(0.2,0.1)$) node {$c_2$};
\draw ($0.5*(X1)+0.5*(Y1)-(0.2,0)$) node {$c$};

\draw ($0.5*(Y1)+0.5*(Z1)+(0.025,0.22)+0.15*(-\horunit,\horunit)$) node {$c_3$};
\draw ($0.5*(Y2)+0.5*(Z2)+(0.025,0.22)+0.15*(-\horunit,\horunit)$) node {$c_1$};

\draw ($0.5*(W1)+0.5*(Z1)+(0.2,-0.1)$) node {$c_2$};
\draw ($0.5*(W2)+0.5*(Z2)+(0.2,0)$) node {$c$};

\foreach \x in {X,W,Y,Z}
\foreach \n in {1,2}
{
\draw [fill] (\x\n) circle [radius=\vxrad];
}
\end{tikzpicture}\hspace{0.6cm}\begin{tikzpicture}\draw [thick,dotted] (0,1.1) -- (0,-1.1);\end{tikzpicture}\hspace{0.6cm}\begin{tikzpicture}[scale=1]
\def\vxrad{0.07cm}
\def\horunit{1.2}
\def\edgelength{0.4}
\def\betweenrows{0.5}

%\draw [white] ($-1*(\horunit,0)$) -- ($41*(\horunit,0)$);

\def\gapratio{1}

\draw [thick,dotted,white] (0,1.1) -- (0,-1.1);

\foreach \num/\parity/\parityy in {1/-1/-1}
{
\coordinate (Y\num) at ($(0,0)+\parityy*\gapratio*(\horunit,0)+(0.5*\parity*\horunit,0.5*\horunit)$);
\coordinate (W\num) at ($(0,0)+\parityy*\gapratio*(\horunit,0)+(-0.5*\parity*\horunit,0.5*\horunit)$);
\coordinate (X\num) at ($(0,0)+\parityy*\gapratio*(\horunit,0)+(0.5*\parity*\horunit,-0.5*\horunit)$);
\coordinate (Z\num) at ($(0,0)+\parityy*\gapratio*(\horunit,0)+(-0.5*\parity*\horunit,-0.5*\horunit)$);
}

%%%%%%%%%%%%%%%NEW LABELS%%%%%%%%%%%%%%%%%%%
\foreach \x/\lab/\offs/\sss in {Z/z/-1.6/5,W/w/1.4/4}
\foreach \num in {1}
{
\draw  ($(\x\num)+\offs*(0,0.25)$) node {$\lab_\num$};
}
\foreach \x/\lab/\offs/\sss in {X/z/-1.6/5,Y/w/1.4/6}
\foreach \num in {2}
{
\draw  ($(\x\num)+\offs*(0,0.25)$) node {$\lab_\num$};
}
%%%%%%%%%%%%%%%NEW LABELS%%%%%%%%%%%%%%%%%%%

\foreach \num/\parity/\parityy in {2/-1/1}
{
\coordinate (Y\num) at ($(0,0)+\parityy*\gapratio*(\horunit,0)+(0.5*\parity*\horunit,0.5*\horunit)$);
%\coordinate (W\num) at ($(0,0)+\parityy*\gapratio*(\horunit,0)+(-0.5*\parity*\horunit,0.5*\horunit)$);
\coordinate (X\num) at ($(0,0)+\parityy*\gapratio*(\horunit,0)+(0.5*\parity*\horunit,-0.5*\horunit)$);
%\coordinate (Z\num) at ($(0,0)+\parityy*\gapratio*(\horunit,0)+(-0.5*\parity*\horunit,-0.5*\horunit)$);
}

\def\drop{0.2};
\draw [rounded corners,thick,black!40,fill=black!20] ($0.5*(Z1)+0.5*(X2)-(0,\drop)$)
-- ($(X2)+(\drop,-\drop)$) -- ($(Y2)+(\drop,\drop)$) -- ($(W1)+(-\drop,\drop)$) -- ($(Z1)+(-\drop,-\drop)$) -- ($0.5*(Z1)+0.5*(X2)-(0,\drop)$);
\draw ($0.25*(W1)+0.25*(Z1)+0.25*(Y2)+0.25*(X2)+(0,0.3)$) node {$V^{\mathrm{abs}}$};

%%%EDGES
\foreach \x/\y/\col/\n in {Z/W/teal/1}%,W/Z/teal/2}
{
%\draw [thick,\col,densely dashed] (\x\n) -- (\y\n);
}

\foreach \x/\y/\col in {X/W/red}
\foreach \n in {1}%,2}
{
\draw [thick,densely dashed] (\x\n) -- (\y\n);
}

\foreach \x/\y/\col in {Y/Z/blue}
\foreach \n in {1}%,2}
{
\draw [thick,densely dotted] (\x\n) -- (\y\n);
}

\foreach \x/\y/\col in {X/Y/magenta}
\foreach \n in {2}
{
\draw [thick] (\x\n) -- (\y\n);
}

%%vertex labels
\foreach \x/\lab/\offs/\sss in {X/x/-1.2/5,Y/y/1/6}%,Z/z/-1.2/7,W/w/1/4}
\foreach \num in {1}
{
\draw  ($(\x\num)+\offs*(0,0.25)$) node {$\lab_\num$};
}
\foreach \x/\lab/\offs/\sss in {X/z/-1.2/5,Y/w/1/6}
\foreach \num in {2}
{
%\draw  ($(\x\num)+\offs*(0,0.25)$) node {$\lab_\num$};
}

%%%colourlabels

%\draw [red] ($0.5*(X2)+0.5*(W2)+(0.025,-0.25)+0.15*(-\horunit,-\horunit)$) node {$c_1$};
\draw ($0.5*(X1)+0.5*(W1)+(0.025,-0.25)+0.15*(-\horunit,-\horunit)$) node {$c_1$};

\draw ($0.5*(X2)+0.5*(Y2)-(0.2,0.1)$) node {$c_2$};
%\draw ($0.5*(X1)+0.5*(Y1)-(0.2,0)$) node {$c_2$};

\draw ($0.5*(Y1)+0.5*(Z1)+(0.025,0.22)+0.15*(-\horunit,\horunit)$) node {$c_3$};
%\draw [blue] ($0.5*(Y2)+0.5*(Z2)+(0.025,0.22)+0.15*(-\horunit,\horunit)$) node {$c_3$};

%\draw ($0.5*(W1)+0.5*(Z1)+(0.2,0)$) node {$c_2$};
%\draw [teal] ($0.5*(W2)+0.5*(Z2)+(0.2,0)$) node {$d$};

\foreach \x in {X,W,Y,Z}
\foreach \n in {1}
{
\draw [fill] (\x\n) circle [radius=\vxrad];
}
\foreach \x in {X,Y}
\foreach \n in {2}
{
\draw [fill] (\x\n) circle [radius=\vxrad];
}
\end{tikzpicture}\hspace{0.6cm}\begin{tikzpicture}\draw [thick,dotted] (0,1.1) -- (0,-1.1);\end{tikzpicture}\hspace{0.6cm}\begin{tikzpicture}[scale=1]
\def\vxrad{0.07cm}
\def\horunit{1.2}
\def\edgelength{0.4}
\def\betweenrows{0.5}

%\draw [white] ($-1*(\horunit,0)$) -- ($41*(\horunit,0)$);

\def\gapratio{1}

\foreach \num/\parity/\parityy in {1/-1/-1}
{
%\coordinate (Y\num) at ($(0,0)+\parityy*\gapratio*(\horunit,0)+(0.5*\parity*\horunit,0.5*\horunit)$);
%\coordinate (W\num) at ($(0,0)+\parityy*\gapratio*(\horunit,0)+(-0.5*\parity*\horunit,0.5*\horunit)$);
\coordinate (X\num) at ($(0,0)+\parityy*\gapratio*(\horunit,0)+(0.5*\parity*\horunit,-0.5*\horunit)$);
\coordinate (Z\num) at ($(0,0)+\parityy*\gapratio*(\horunit,0)+(-0.5*\parity*\horunit,-0.5*\horunit)$);
}
\foreach \num/\parity/\parityy in {1/-1/-1,2/-1/1}
{
\coordinate (Y\num) at ($(0,0)+\parityy*\gapratio*(\horunit,0)+(0.5*\parity*\horunit,0.5*\horunit)$);
\coordinate (W\num) at ($(0,0)+\parityy*\gapratio*(\horunit,0)+(-0.5*\parity*\horunit,0.5*\horunit)$);
\coordinate (X\num) at ($(0,0)+\parityy*\gapratio*(\horunit,0)+(0.5*\parity*\horunit,-0.5*\horunit)$);
\coordinate (Z\num) at ($(0,0)+\parityy*\gapratio*(\horunit,0)+(-0.5*\parity*\horunit,-0.5*\horunit)$);
}

\def\drop{0.2};
\draw [rounded corners,thick,black!40,fill=black!20] ($0.5*(Z1)+0.5*(X2)-(0,\drop)$)
-- ($(X2)+(\drop,-\drop)$) -- ($(Y2)+(\drop,\drop)$) -- ($(W1)+(-\drop,\drop)$) -- ($(Z1)+(-\drop,-\drop)$) -- ($0.5*(Z1)+0.5*(X2)-(0,\drop)$);
\draw ($0.25*(W1)+0.25*(Z1)+0.25*(Y2)+0.25*(X2)+(0,0.3)$) node {$V^{\mathrm{abs}}$};

%%%EDGES
\foreach \x/\y/\col/\n in {Z/W/teal/1}
{
\draw [thick] (\x\n) -- (\y\n);
}

\foreach \x/\y/\col in {X/W/red}
\foreach \n in {2}
{
\draw [thick,densely dotted] (\x\n) -- (\y\n);
}

%%%%%%%%%%%%%%%NEW LABELS%%%%%%%%%%%%%%%%%%%
\foreach \x/\lab/\offs/\sss in {Z/z/-1.6/5,W/w/1.4/4}
\foreach \num in {1}
{
\draw  ($(\x\num)+\offs*(0,0.25)$) node {$\lab_\num$};
}
\foreach \x/\lab/\offs/\sss in {X/z/-1.6/5,Y/w/1.4/6}
\foreach \num in {2}
{
\draw  ($(\x\num)+\offs*(0,0.25)$) node {$\lab_\num$};
}
%%%%%%%%%%%%%%%NEW LABELS%%%%%%%%%%%%%%%%%%%

\foreach \x/\y/\col in {Y/Z/blue}
\foreach \n in {2}
{
\draw [thick,densely dashed] (\x\n) -- (\y\n);
}

\foreach \x/\y/\col in {X/Y/magenta}
\foreach \n in {2}
{
%\draw [thick,\col] (\x\n) -- (\y\n);
}

%%vertex labels
\foreach \x/\lab/\offs/\sss in {Z/z/-1.2/7,W/w/1/4}%X/x/-1.2/5,Y/y/1/6,
\foreach \num in {1}
{
%\draw  ($(\x\num)+\offs*(0,0.25)$) node {$\lab_\num$};
}
\foreach \x/\lab/\offs/\sss in {Z/x/-1.2/7,W/y/1/4}%X/z/-1.2/5,Y/w/1/6,
\foreach \num in {2}
{
\draw  ($(\x\num)+\offs*(0,0.25)$) node {$\lab_\num$};
}

%%%colourlabels

\draw ($0.5*(X2)+0.5*(W2)+(0.025,-0.25)+0.15*(-\horunit,-\horunit)$) node {$c_3$};
%\draw [red] ($0.5*(X1)+0.5*(W1)+(0.025,-0.25)+0.15*(-\horunit,-\horunit)$) node {$c_1$};

%\draw [magenta] ($0.5*(X2)+0.5*(Y2)-(0.2,0)$) node {$c_3$};
%\draw [magenta] ($0.5*(X1)+0.5*(Y1)-(0.2,0)$) node {$c_2$};

%\draw [blue] ($0.5*(Y1)+0.5*(Z1)+(0.025,0.22)+0.15*(-\horunit,\horunit)$) node {$c_3$};
\draw ($0.5*(Y2)+0.5*(Z2)+(0.025,0.22)+0.15*(-\horunit,\horunit)$) node {$c_1$};

\draw  ($0.5*(W1)+0.5*(Z1)+(0.2,-0.1)$) node {$c_2$};
%\draw [teal] ($0.5*(W2)+0.5*(Z2)+(0.2,0)$) node {$d$};

\foreach \x in {X,W,Y,Z}
\foreach \n in {2}
{
\draw [fill] (\x\n) circle [radius=\vxrad];
}
\foreach \x in {W,Z}%X,Y
\foreach \n in {1}
{
\draw [fill] (\x\n) circle [radius=\vxrad];
}
\end{tikzpicture}
\end{center}

\vspace{-0.3cm}

\caption{On the left a structure allowing $(V^{\mathrm{abs}},C^{\mathrm{abs}})$ to absorb either $\{x_1,y_1\}$ or $\{x_2,y_2\}$, where $C^{\mathrm{abs}}=\{c_1,c_2,c_3\})$. That is, there is both a matching with vertex set $V^{\mathrm{abs}}\cup \{x_1,y_1\}$ and colour set $C^{\mathrm{abs}}$ (in the middle) and $V^{\mathrm{abs}}\cup \{x_2,y_2\}$ and colour set $C^{\mathrm{abs}}$ (on the right). In a colouring arising from a group, this is possible as $x_1y_1$ and $x_2y_2$ have the same colour.}\label{fig:switcher}
\end{figure}

This is accomplished by M\"uyesser and Pokrovskiy in~\cite{muyesser2022random} with $\ell_0=o(n)$ and $\ell_1=n^{1-\varepsilon}$ for some small fixed $\varepsilon>0$. (In fact, they prove something stronger that we have weakened to cohere with our subsequent discussion.)
To build this absorber, they use distributive absorption (as mentioned in Section~\ref{sec:decomp}), an efficient method to build a global absorption property from small local absorbers. To avoid excessive detail here, we will comment no further than to say that the starting point is usually to robustly find small absorbers. For example, here, it would be useful to find a small absorber that can absorb either $\{x_1,y_1\}$ or $\{x_2,y_2\}$ (two balanced sets). The discussion around \eqref{eqn:recall} implies that this is only possible if $x_1+y_1=x_2+y_2$, or, equally, if $x_1y_1$ and $x_2y_2$ have the same colour in $H=H(G)$. (Note we are giving a minimal useful example as if we could absorb $\{x_1\}$ and $\{x_2\}$ then in fact $x_1=x_2$.) However, subject only to the condition that $x_1+y_1=x_2+y_2$ ($=c$, say) we can find such a small absorber in many ways (as depicted in Figure~\ref{fig:switcher}). Picking any two colours $c_1$ and $c_2$ (avoiding $O(1)$ bad options so the following is possible), label vertices $w_1,z_1$ and a colour $c_3$ such that $x_1w_1z_1y_1$ is a path in $H$ with colours $c_1,c_2$ and $c_3$ in that order, and $c,c_1,c_2,c_3$ are all distinct. Label vertices $w_2,z_2$ such that  $x_2w_2z_2y_2$ is a path in $H$ with colours $c_1,c_2$ and $c_3$ in that order, noting that this is possible as $H=H(G)$ and $x_2y_2$ is the same colour as $x_1y_1$. Then, as long as we chose $c_1,c_2$ so that these paths are vertex-disjoint,  it is easy to see that, setting $V^{\mathrm{abs}}=\{w_1,z_1,w_2,z_2\}$ and $C^{\mathrm{abs}}=\{c_1,c_2,c_3\}$, $(V^{\mathrm{abs}},C^{\mathrm{abs}})$ can absorb either $\{x_1,y_1\}$ or $\{x_2,y_2\}$ (see Figure~\ref{fig:switcher}). %(Note that technically $(\emptyset,\{c\})$ absorbs either of these sets, but this is not a robust constructiopn )

There are many complications dealt with in~\cite{muyesser2022random}, not least of all of course when $G$ is non-abelian, but we hope this gives some indication of the starting point from which a global absorber can be created. We return to this in Section~\ref{sec:partial} to discuss the challenges in attempting this more generally for colourings not generated from a group.  Constructing absorbers for zero-sum sets in this manner is a natural extension in the development of absorption techniques, and appears also in the recent work by Bowtell and Keevash~\cite{bowtell2021n} on the famous $n$-queens problem.

In Section~\ref{sec:decomp}, we described work by Gould and Kelly~\cite{gould2023hamilton} finding Hamilton transversals in a random Latin square with high probability. There, a Hamilton transversal corresponded to a full transversal where the natural permutation it defines is a cyclic permutation.
More generally, we could ask for a full transversal in a random Latin square of order $n$ where this permutation has any specific cycle type corresponding to a permutation of an $n$-element set. If there are few fixed points then it seems plausible this can be done with high probability.
In the non-random case, M\"uyesser~\cite{muyesser2023cycle} has conjectured that if $k\geq 3$ is such that $k|(n-1)$, then for any group $G$ of order $n$ satisfying the Hall-Paige condition the Latin square $L(G)$ has a full transversal whose permutation corresponds to a single fixed point and otherwise has cycles of length $k$. (The case $k=n-1$ can be seen to correspond to a Hamilton transversal.)
In related work, M\"uyesser~\cite{muyesser2023cycle} proved for large groups a conjecture of Friedlander, Gordon and Tannenbaum~\cite{friedlander1978group} from 1981 that, for $k\geq 2$ dividing $(n-1)$, if an \emph{abelian} group $G$ of order $n$ satisfies the Hall-Paige condition then the Latin square formed like $L(G)$ but with entry $u-v$ in the cell indexed by row $u$ and column $v$, for each $u,v\in G$, has a transversal whose permutation corresponds to a single fixed point and otherwise has cycles of length $k$.

%, who conjectured that it should hold more widely for any group satisfying the Hall-Paige conjecture as long as $k\geq 3$.

%in 1981 Friedlander, Gordon and Tannenbaum~\cite{friedlander1978group} conjectured that if an \emph{abelian} group $G$ of order $n$ satisfies the Hall-Paige condition, and $k\geq 2$ satisfies $k|(n-1)$, then $L(G)$ has a full transversal whose permutation corresponds to a single fixed point and otherwise has cycles of length $k$. (The case $k=n-1$ can be seen to correspond to a Hamilton transversal.) This conjecture has recently been proved for large groups by M\"uyesser~\cite{muyesser2023cycle}, who conjectured that it should hold more widely for any group satisfying the Hall-Paige conjecture as long as $k\geq 3$.

%%%%%%%%%%%%%%%%%%%%%%%%%%%%%%%%%%%%%%%%%%%%%%%%%%%%%%%%%%%%%%%%%%%%%%%%%%%%%%%%%%%%%%%%%%%%%%%%%%%%%%%%%%%%%%%%%%%%%%%%%%
%%%%%%%%%%%%%%%%%%%%%%%%%%%%%%%%%%%%%%%%%%%%%%%%%%%%%%%%%%%%%%%%%%%%%%%%%%%%%%%%%%%%%%%%%%%%%%%%%%%%%%%%%%%%%%%%%%%%%%%%%%
%%%%%%%%%%%%%%%%%%%%%%%%%%%%%%%%%%%%%%%%%%%%%%%%%%%%%%%%%%%%%%%%%%%%%%%%%%%%%%%%%%%%%%%%%%%%%%%%%%%%%%%%%%%%%%%%%%%%%%%%%%
%%%%%%%%%%%%%%%%%%%%%%%%%%%%%%%%%%%%%%%%%%%%%%%%%%%%%%%%%%%%%%%%%%%%%%%%%%%%%%%%%%%%%%%%%%%%%%%%%%%%%%%%%%%%%%%%%%%%%%%%%%

\section{Partial transversals in Latin squares}\label{sec:partial}
Shortly after the formation of the Ryser-Brualdi-Stein conjecture (or, more precisely, the portion due to Ryser~\cite{Ryser}), Koksma~\cite{koksma} gave a simple argument that every Latin square of order $n$ has a partial transversal with at least $2n/3$ cells. Drake~\cite{drake} improved this bound to $3n/4$, before an approximate form of the conjecture was proved in 1978 independently by
Brouwer, De Vries and Wieringa~\cite{brouwer1978lower} and  Woolbright~\cite{woolbright}, each showing that every Latin square of order $n$ has a partial transversal with at least $n-\sqrt{n}$ cells. In 1982, Shor~\cite{shor} gave a proof improving this bound to $n-O(\log^2n)$, though this paper contained an error that was only noticed and fixed in 2008 by Hatami and Shor~\cite{hatamishor}, using effectively the original approach.
Recently, Keevash, Pokrovskiy, Sudakov and Yepremyan~\cite{KPSY} improved this bound, which had essentially stood for 40 years, to show any Latin square of order $n$ has a partial transversal with $n-O(\log n/\log\log n)$ cells.
Very recently, Montgomery~\cite{montgomeryrbs} showed that every sufficiently large Latin square of order $n$ has a partial transversal with $n-1$ cells. While this confirms the Ryser-Brualdi-Stein for large even $n$, it seems very difficult to find a full transversal in a Latin square with large odd order using these methods and certainly further ideas would be needed.

To discuss the proof of the bounds in~\cite{KPSY,montgomeryrbs}, we will use the rainbow subgraph formulation of the problem. Suppose, then, that $G$ is a bipartite copy of $K_{n,n}$ which is properly coloured with $n$ colours. To find a find a large rainbow matching in $G$, Keevash, Pokrovskiy, Sudakov and Yepremyan~\cite{KPSY} essentially began with a large random almost-perfect rainbow matching, before successively modifying it, each time getting a rainbow matching that is one edge larger, until it contains $n-O(\log n/\log\log n)$ edges. Suppose the initial large rainbow matching is $M$. Suppose further that there are vertices $x,y\in V(G)\setminus V(M)$ which are in different vertex classes in $G$, and an $x,y$-path $P$ in $G$ such that the even edges in this path are in $M$ and the odd edges (starting with the edge containing $x$) all have different colours not used on $M$. Then, removing the even edges of $P$ from $M$ and adding the odd edges gives a rainbow matching in $G$ with one more edge. In~\cite{KPSY}, the initial large rainbow rainbow matching $M$ is found using the R\"odl nibble (as discussed briefly in Section~\ref{sec:decomp}), and has $n-n^{1-\alpha}$ edges for some fixed small $\alpha>0$. The randomness of the initial matching $M$ implies (via some expansion properties) essentially that for every pair of vertices not in the matching and every set of $K\log n/\log\log n$ colours (for some fixed constant $K$) not on $M$ such a path $P$ can be found. Moreover, this property is sufficiently robust that it can be used to make iterative adjustments (with some additional care) until the rainbow matching uses all but $K\log n/\log\log n$ of the colours of $G$.

This brief outline belies the challenges overcome in~\cite{KPSY}, particularly in achieving a bound that seems likely to be optimal using these or related methods. Indeed, if $d=\log n/2\log\log n$ and $C$ is a set of $d$ colours, then however a matching $M$ is chosen in $G$ without using the colours in $C$, for any vertex $v\in V(G)$, there are at most $2d^d=o(\sqrt{n})$ paths from $v$ which alternate between edges with colour in $C$ and edges in $M$ in the manner of the path $P$ described above. Thus, such paths certainly cannot exist between an arbitrary pair of vertices. More broadly, the bound $n-O(\log n/\log\log n)$ is plausibly a natural barrier for any method that does not alter its approach based on the specific colour of $G$, and the underlying algebraic properties it might have (see, for example, the examples described in Section~\ref{sec:examples}).

A significant part of the work in~\cite{montgomeryrbs}, then, is to find some way to identify and exploit (approximate) algebraic properties in the colouring of $G$. Following this, the approach taken is an implementation of the absorption method, where, as in the outline of the work in~\cite{muyesser2023cycle} by M\"uyesser and Pokrovskiy in Section~\ref{sec:full}, an absorption structure is found and set aside, before a large rainbow matching is found disjointly. Finding a large rainbow matching disjoint from the absorption structure in~\cite{montgomeryrbs} is possible if a substantial set of random colours and vertices is reserved for this, following work by Montgomery, Pokrovskiy and Sudakov~\cite{montgomery2021proof}.

The key, then, is in developing some absorption structure. We recall from Section~\ref{sec:full} that we cannot hope to find an absorption structure that can absorb any small set of vertices. In contrast to the previous methods involving the absorption of `zero-sum' sets (i.e., those in~\cite{bowtell2021n,muyesser2022random}), here we do not have a structure with a well-determined, known, and exact algebraic structure.
Some condition is thus needed to define the sets we will be able to absorb. In~\cite{montgomeryrbs}, a colour $c_0$ is chosen (fairly arbitrarily) to function as an `identity colour' and the condition for absorption is that the vertex set must be the vertex set of a matching of colour-$c_0$ edges. (A few edges are deleted for this to be true, but here we gloss over this for simplicity.) More precisely,  for some $\ell_0,\ell_1\leq n$ we look for a balanced set $V^{\mathrm{abs}}$ of $2\ell_0$ vertices and a set $C^{\mathrm{abs}}$ of $\ell_0+\ell_1$
colours such that, for any balanced set $W\subset V(H)\setminus V^{\mathrm{abs}}$  \emph{which is the vertex set of a matching of $\ell_1$ edges with colour $c_0$}, there is a rainbow matching in $G$ with vertex set $V^{\mathrm{abs}}\cup W$ and colour set $C^{\mathrm{abs}}$.

\begin{figure}
\begin{center}
\begin{tikzpicture}[scale=1]
\def\vxrad{0.07cm}
\def\horunit{1.2}
\def\edgelength{0.4}
\def\betweenrows{0.5}

%\draw [white] ($-1*(\horunit,0)$) -- ($41*(\horunit,0)$);

\def\gapratio{1}

\foreach \num/\parity/\parityy in {1/-1/-1,2/-1/1}
{
\coordinate (Y\num) at ($(0,0)+\parityy*\gapratio*(\horunit,0)+(0.5*\parity*\horunit,0.5*\horunit)$);
\coordinate (W\num) at ($(0,0)+\parityy*\gapratio*(\horunit,0)+(-0.5*\parity*\horunit,0.5*\horunit)$);
\coordinate (X\num) at ($(0,0)+\parityy*\gapratio*(\horunit,0)+(0.5*\parity*\horunit,-0.5*\horunit)$);
\coordinate (Z\num) at ($(0,0)+\parityy*\gapratio*(\horunit,0)+(-0.5*\parity*\horunit,-0.5*\horunit)$);
}

\def\drop{0.2};
\draw [rounded corners,thick,black!40,fill=black!10] ($0.5*(Z1)+0.5*(X2)-(0,\drop)$)
-- ($(X2)+(\drop,-\drop)$) -- ($(Y2)+(\drop,\drop)$) -- ($(W1)+(-\drop,\drop)$) -- ($(Z1)+(-\drop,-\drop)$) -- ($0.5*(Z1)+0.5*(X2)-(0,\drop)$);
\draw ($0.25*(W1)+0.25*(Z1)+0.25*(Y2)+0.25*(X2)+(0,0.3)$) node {$V^{\mathrm{abs}}$};

%%%EDGES
\foreach \x/\y/\col/\n in {Z/W/teal/1}
{
\draw [thick] (\x\n) -- (\y\n);
}
\foreach \x/\y/\col/\n in {X/Y/teal/2}
{
%\draw [thick, dash dot] (\x\n) -- (\y\n);
\draw [thick] (\x\n) -- ($0.6*(\x\n)+0.4*(\y\n)$);
\draw [thick] ($0.55*(\x\n)+0.45*(\y\n)$) -- ($0.45*(\x\n)+0.55*(\y\n)$);
\draw [thick] ($0.4*(\x\n)+0.6*(\y\n)$) -- (\y\n);
}

\foreach \x/\y/\col in {X/W/red}
\foreach \n in {1}
{
\draw [thick,densely dashed] (\x\n) -- (\y\n);
}

\foreach \x/\y/\col in {Y/Z/red}
\foreach \n in {2}
{
\draw [thick,densely dashed] (\x\n) -- (\y\n);
}

\foreach \x/\y/\col in {Y/Z/blue}
\foreach \n in {1}
{
\draw [thick,densely dotted] (\x\n) -- (\y\n);
}
\foreach \x/\y/\col in {X/W/blue}
\foreach \n in {2}
{
\draw [thick,dotted] (\x\n) -- (\y\n);
}

\foreach \x/\y/\col in {X/Y/magenta}
\foreach \n in {1}
{
\draw [thick,black!50] (\x\n) -- (\y\n);
}
\foreach \x/\y/\col in {W/Z/magenta}
\foreach \n in {2}
{
\draw [thick,black!50] (\x\n) -- (\y\n);
}

%%vertex labels
\foreach \x/\lab/\offs/\sss in {X/x/-1.2/5,Y/y/1/6}%,Z/z/-1.2/7,W/w/1/4}
\foreach \num in {1}
{
\draw  ($(\x\num)+\offs*(0,0.25)$) node {$\lab_\num$};
}
\foreach \x/\lab/\offs/\sss in {Z/x/-1.2/7,W/y/1/4}%X/z/-1.2/5,Y/w/1/6,
\foreach \num in {2}
{
\draw  ($(\x\num)+\offs*(0,0.25)$) node {$\lab_\num$};
}

%%%colourlabels

\draw ($0.5*(X2)+0.5*(W2)+(0.025,-0.25)+0.15*(-\horunit,-\horunit)$) node {$c_3'$};
\draw ($0.5*(X1)+0.5*(W1)+(0.025,-0.25)+0.15*(-\horunit,-\horunit)$) node {$c_1$};

\draw ($0.5*(X2)+0.5*(Y2)-(0.2,0.1)$) node {$c_2$};
\draw ($0.5*(X1)+0.5*(Y1)-(0.2,0)$) node {$c$};

\draw ($0.5*(Y1)+0.5*(Z1)+(0.025,0.22)+0.15*(-\horunit,\horunit)$) node {$c_3$};
\draw ($0.5*(Y2)+0.5*(Z2)+(0.025,0.22)+0.15*(-\horunit,\horunit)$) node {$c_1$};

\draw ($0.5*(W1)+0.5*(Z1)+(0.2,-0.1)$) node {$c_2$};
\draw ($0.5*(W2)+0.5*(Z2)+(0.2,0)$) node {$c$};

\foreach \x in {X,W,Y,Z}
\foreach \n in {1,2}
{
\draw [fill] (\x\n) circle [radius=\vxrad];
}
\end{tikzpicture}\hspace{0cm}\begin{tikzpicture}\draw [thick,dotted] (0,1.1) -- (0,-1.1);\end{tikzpicture}\hspace{0.3cm}\begin{tikzpicture}[scale=1]
\def\vxrad{0.07cm}
\def\horunit{1.2}
\def\edgelength{0.4}
\def\betweenrows{0.5}

%\draw [white] ($-1*(\horunit,0)$) -- ($41*(\horunit,0)$);

\def\gapratio{1}

\draw [thick,dotted,white] (0,1.1) -- (0,-1.1);

\foreach \num/\parity/\parityy in {1/-1/-1,2/-1/1}
{
\coordinate (Y\num) at ($(0,0)+\parityy*\gapratio*(\horunit,0)+(0.5*\parity*\horunit,0.5*\horunit)$);
\coordinate (W\num) at ($(0,0)+\parityy*\gapratio*(\horunit,0)+(-0.5*\parity*\horunit,0.5*\horunit)$);
\coordinate (X\num) at ($(0,0)+\parityy*\gapratio*(\horunit,0)+(0.5*\parity*\horunit,-0.5*\horunit)$);
\coordinate (Z\num) at ($(0,0)+\parityy*\gapratio*(\horunit,0)+(-0.5*\parity*\horunit,-0.5*\horunit)$);
}

\def\drop{0.2};
\def\slidd{1};
\draw [rounded corners,thick,black!40,fill=black!10] ($0.5*(X1)+0.5*(Z2)-(0,\drop)$)
-- ($(Z2)+(\slidd*\drop,-\drop)$) -- ($(W2)+(\slidd*\drop,\drop)$) -- ($(Y1)+(-\slidd*\drop,\drop)$) -- ($(X1)+(-\slidd*\drop,-\drop)$) -- ($0.5*(X1)+0.5*(Z2)-(0,\drop)$);
%\draw ($0.25*(W1)+0.25*(Z1)+0.25*(Y2)+0.25*(X2)+(0,0.3)$) node {$V$};

%%%EDGES
\foreach \x/\y/\col/\n in {Z/W/teal/1}
{
%\draw [thick] (\x\n) -- (\y\n);
}
\foreach \x/\y/\col/\n in {X/Y/teal/2}
{
%\draw [thick, dash dot] (\x\n) -- (\y\n);
\draw [thick] (\x\n) -- ($0.6*(\x\n)+0.4*(\y\n)$);
\draw [thick] ($0.55*(\x\n)+0.45*(\y\n)$) -- ($0.45*(\x\n)+0.55*(\y\n)$);
\draw [thick] ($0.4*(\x\n)+0.6*(\y\n)$) -- (\y\n);
}

\foreach \x/\y/\col in {X/W/red}
\foreach \n in {1}
{
\draw [thick,densely dashed] (\x\n) -- (\y\n);
}

\foreach \x/\y/\col in {Y/Z/red}
\foreach \n in {2}
{
%\draw [thick,densely dashed] (\x\n) -- (\y\n);
}

\foreach \x/\y/\col in {Y/Z/blue}
\foreach \n in {1}
{
\draw [thick,densely dotted] (\x\n) -- (\y\n);
}
\foreach \x/\y/\col in {X/W/blue}
\foreach \n in {2}
{
%\draw [thick,densely dotted] (\x\n) -- (\y\n);
}

\foreach \x/\y/\col in {X/Y/magenta}
\foreach \n in {1}
{
%\draw [thick,black!50] (\x\n) -- (\y\n);
}
\foreach \x/\y/\col in {W/Z/magenta}
\foreach \n in {2}
{
\draw [thick,black!50] (\x\n) -- (\y\n);
}

%%vertex labels
\foreach \x/\lab/\offs/\sss in {X/x/-1.2/5,Y/y/1/6}%,Z/z/-1.2/7,W/w/1/4}
\foreach \num in {1}
{
%\draw  ($(\x\num)+\offs*(0,0.25)$) node {$\lab_\num$};
}
\foreach \x/\lab/\offs/\sss in {Z/x/-1.2/7,W/y/1/4}%X/z/-1.2/5,Y/w/1/6,
\foreach \num in {2}
{
%\draw  ($(\x\num)+\offs*(0,0.25)$) node {$\lab_\num$};
}

%%%colourlabels

%\draw ($0.5*(X2)+0.5*(W2)+(0.025,-0.25)+0.15*(-\horunit,-\horunit)$) node {$c_3$};
\draw ($0.5*(X1)+0.5*(W1)+(0.025,-0.25)+0.15*(-\horunit,-\horunit)$) node {$c_1$};

\draw ($0.5*(X2)+0.5*(Y2)-(0.2,0.1)$) node {$c_2$};
%\draw ($0.5*(X1)+0.5*(Y1)+(0.2,0)$) node {$c$};

\draw ($0.5*(Y1)+0.5*(Z1)+(0.025,0.22)+0.15*(-\horunit,\horunit)$) node {$c_3$};
%\draw ($0.5*(Y2)+0.5*(Z2)+(0.025,0.22)+0.15*(-\horunit,\horunit)$) node {$c_1$};

%\draw ($0.5*(W1)+0.5*(Z1)+(0.2,-0.1)$) node {$c_2$};
\draw ($0.5*(W2)+0.5*(Z2)+(-0.2,0)$) node {$c$};

\draw ($0.25*(W1)+0.25*(Z1)+0.25*(Y2)+0.25*(X2)-(0.2,0)$) node {$\hat{V}^{\mathrm{abs}}$};

\foreach \x in {X,W,Y,Z}
\foreach \n in {1,2}
{
\draw [fill] (\x\n) circle [radius=\vxrad];
}
\end{tikzpicture}\hspace{0.3cm}\begin{tikzpicture}\draw [thick,dotted] (0,1.1) -- (0,-1.1);\end{tikzpicture}\hspace{0.3cm}\begin{tikzpicture}[scale=1]
\def\vxrad{0.07cm}
\def\horunit{1.2}
\def\edgelength{0.4}
\def\betweenrows{0.5}

%\draw [white] ($-1*(\horunit,0)$) -- ($41*(\horunit,0)$);

\def\gapratio{1}

\draw [thick,dotted,white] (0,1.1) -- (0,-1.1);

\foreach \num/\parity/\parityy in {1/-1/-1,2/-1/1}
{
\coordinate (Y\num) at ($(0,0)+\parityy*\gapratio*(\horunit,0)+(0.5*\parity*\horunit,0.5*\horunit)$);
\coordinate (W\num) at ($(0,0)+\parityy*\gapratio*(\horunit,0)+(-0.5*\parity*\horunit,0.5*\horunit)$);
\coordinate (X\num) at ($(0,0)+\parityy*\gapratio*(\horunit,0)+(0.5*\parity*\horunit,-0.5*\horunit)$);
\coordinate (Z\num) at ($(0,0)+\parityy*\gapratio*(\horunit,0)+(-0.5*\parity*\horunit,-0.5*\horunit)$);
}

\def\drop{0.2};
\def\slidd{1};
\draw [rounded corners,thick,black!40,fill=black!10] ($0.5*(X1)+0.5*(Z2)-(0,\drop)$)
-- ($(Z2)+(\slidd*\drop,-\drop)$) -- ($(W2)+(\slidd*\drop,\drop)$) -- ($(Y1)+(-\slidd*\drop,\drop)$) -- ($(X1)+(-\slidd*\drop,-\drop)$) -- ($0.5*(X1)+0.5*(Z2)-(0,\drop)$);
%\draw ($0.25*(W1)+0.25*(Z1)+0.25*(Y2)+0.25*(X2)+(0,0.3)$) node {$V$};

%%%EDGES
\foreach \x/\y/\col/\n in {Z/W/teal/1}
{
\draw [thick] (\x\n) -- (\y\n);
}
\foreach \x/\y/\col/\n in {X/Y/teal/2}
{
%\draw [thick, dash dot] (\x\n) -- (\y\n);
%\draw [thick] (\x\n) -- ($0.6*(\x\n)+0.4*(\y\n)$);
%\draw [thick] ($0.55*(\x\n)+0.45*(\y\n)$) -- ($0.45*(\x\n)+0.55*(\y\n)$);
%\draw [thick] ($0.4*(\x\n)+0.6*(\y\n)$) -- (\y\n);
}

\foreach \x/\y/\col in {X/W/red}
\foreach \n in {1}
{
%\draw [thick,densely dashed] (\x\n) -- (\y\n);
}

\foreach \x/\y/\col in {Y/Z/red}
\foreach \n in {2}
{
\draw [thick,densely dashed] (\x\n) -- (\y\n);
}

\foreach \x/\y/\col in {Y/Z/blue}
\foreach \n in {1}
{
%\draw [thick,densely dotted] (\x\n) -- (\y\n);
}
\foreach \x/\y/\col in {X/W/blue}
\foreach \n in {2}
{
\draw [thick,dotted] (\x\n) -- (\y\n);
}

\foreach \x/\y/\col in {X/Y/magenta}
\foreach \n in {1}
{
\draw [thick,black!50] (\x\n) -- (\y\n);
}
\foreach \x/\y/\col in {W/Z/magenta}
\foreach \n in {2}
{
%\draw [thick,black!50] (\x\n) -- (\y\n);
}

%%vertex labels
\foreach \x/\lab/\offs/\sss in {X/x/-1.2/5,Y/y/1/6}%,Z/z/-1.2/7,W/w/1/4}
\foreach \num in {1}
{
%\draw  ($(\x\num)+\offs*(0,0.25)$) node {$\lab_\num$};
}
\foreach \x/\lab/\offs/\sss in {Z/x/-1.2/7,W/y/1/4}%X/z/-1.2/5,Y/w/1/6,
\foreach \num in {2}
{
%\draw  ($(\x\num)+\offs*(0,0.25)$) node {$\lab_\num$};
}

%%%colourlabels

\draw ($0.5*(X2)+0.5*(W2)+(0.025,-0.25)+0.15*(-\horunit,-\horunit)$) node {$c_3'$};
%\draw ($0.5*(X1)+0.5*(W1)+(0.025,-0.25)+0.15*(-\horunit,-\horunit)$) node {$c_1$};

%\draw ($0.5*(X2)+0.5*(Y2)-(0.2,0.1)$) node {$c_2'$};
\draw ($0.5*(X1)+0.5*(Y1)+(0.2,0)$) node {$c$};

%\draw ($0.5*(Y1)+0.5*(Z1)+(0.025,0.22)+0.15*(-\horunit,\horunit)$) node {$c_3$};
\draw ($0.5*(Y2)+0.5*(Z2)+(0.025,0.22)+0.15*(-\horunit,\horunit)$) node {$c_1$};

\draw ($0.5*(W1)+0.5*(Z1)+(0.2,-0.1)$) node {$c_2$};
%\draw ($0.5*(W2)+0.5*(Z2)+(0.2,0)$) node {$c$};

\draw ($0.25*(W1)+0.25*(Z1)+0.25*(Y2)+0.25*(X2)+(0.2,0)$) node {$\hat{V}^{\mathrm{abs}}$};

\foreach \x in {X,W,Y,Z}
\foreach \n in {1,2}
{
\draw [fill] (\x\n) circle [radius=\vxrad];
}
\end{tikzpicture}
\end{center}
\caption{On the left, in contrast to Figure~\ref{fig:switcher}, with $C^{\mathrm{abs}}=\{c_1,c_2,c_3\}$, $({V}^{\mathrm{abs}},C^{\mathrm{abs}})$ can absorb $\{x_1,y_1\}$ but \emph{not} $\{x_2,y_2\}$ if $c_3'\neq c_3$.
Instead, with $\hat{V}^{\mathrm{abs}}=V^{\mathrm{abs}}\cup\{x_1,y_1,x_2,y_2\}$ and $\hat{C}^{\mathrm{abs}}=\{c,c_1,c_2\}$, $(\hat{V}^{\mathrm{abs}},\hat{C}^{\mathrm{abs}})$ can absorb either $c_3$ or $c_3'$, as depicted by the matchings in the middle and on the right, whose colour sets differ on $\{c_3,c_3'\}$.}\label{fig:colswitcher}
\end{figure}
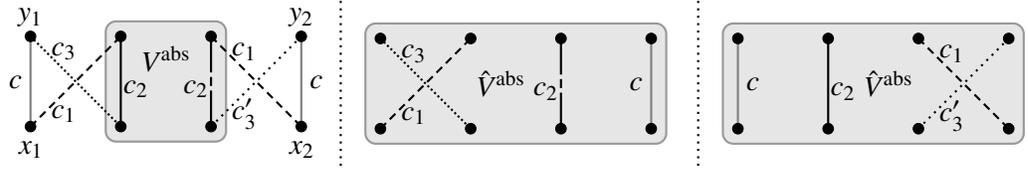

In~\cite{montgomeryrbs}, this absorption structure is found using distributive absorption, so, as in our discussion in Section~\ref{sec:full}, a good starting point would be to find an absorber capable of absorbing one of two sets of two vertices with `zero sum' -- here, instead, the `zero sum' condition is replaced by the set of two vertices being the vertex set of an edge with the identity colour, $c_0$.
However, trying the construction in Section~\ref{sec:full} with $c=c_0$ (see Figure~\ref{fig:switcher}) to find an absorber here does not work as, without the algebraic structure conferred by a group, we may not have that the third edge of the $x_2,y_2$-path has colour $c_3$, but, perhaps, some other colour $c_3'$ (see Figure~\ref{fig:colswitcher}).
Crucially, however, we can instead choose from what we have found a structure that can absorb either the colour $c_3$ or the colour $c_3'$ (as depicted in Figure~\ref{fig:colswitcher}). That is, using the labelling in Figures~\ref{fig:switcher} and~\ref{fig:colswitcher}, taking $\hat{V}^{\mathrm{abs}}$
to be the set of all the vertices used and $\hat{C}^{\mathrm{abs}}$ to be the set of all the colours used except for $c_3$ and $c_3'$, we have that there is a rainbow matching with vertex set $\hat{V}^{\mathrm{abs}}$ and colour set $\hat{C}^{\mathrm{abs}}\cup\{c_3\}$ and a rainbow matching with vertex set $\hat{V}^{\mathrm{abs}}$ and colour set $\hat{C}^{\mathrm{abs}}\cup\{c_3'\}$ (see Figure~\ref{fig:colswitcher}).
Note that if we knew $G$ has many such small absorbers capable of absorbing either $c_3$ or $c_3'$, then if we tried the previous absorber construction for $\{x_1,y_1\}$ and $\{x_2,y_2\}$ and encountered the same problem we could now solve it. Indeed, if we take, using new colours and vertices, a small absorber that can absorb either $c_3$ or $c_3'$, then, in combination with the construction that failed because $c_3\neq c_3'$ (adding both $c_3$ and $c_3'$), it is relatively easy to see that this can absorb either $\{x_1,y_1\}$ or $\{x_2,y_2\}$.

This is a small toe-hold into creating an absorption structure, but a crucial one, and it motivates how we might study the colouring of $G$ to gather something representative of some of its (possible) algebraic structure. Looking at the pairs of colours, we consider how many different small absorbers can absorb either one (as in Figure~\ref{fig:colswitcher}). A good way to consider this to take a complete auxiliary graph $K$ with the colour set of $G$ as its vertex set where each edge is weighted by the number of absorbers in $G$ with the construction in Figure~\ref{fig:colswitcher} which can absorb either one of those two colours.
If two colours $c$ and $d$ are connected by a short path of high weight edges in $K$, then, loosely speaking, we can chain disjoint absorbers together along this path to create an absorber that can absorb either $c$ or $d$. Ideally, we want to know that, the majority of the time, when we attempt the construction in Figure~\ref{fig:colswitcher}, if it fails because $c'_3\neq c_3$, then we will be able to find many small absorbers that can absorb $c'_3$ or $c_3$. To do this, in~\cite{montgomeryrbs}, most of the weight of $K$ is essentially covered by well-connected subgraphs whose edges have roughly equal weight in such a way that no colour appears in too many subgraphs. This is done using techniques involving \emph{sublinear expansion}, which allows well-connected graphs to be found where the length of the connections can be controlled (so that any absorbers created using such a connection in $K$ are not too large). Sublinear expansion, as introduced by Koml\'os and Szemer\'edi~\cite{K-Sz-1,K-Sz-2}, is used here as the non-zero-weighted edges of $K$ might be sparse, yet carry a meaningful weight. For more on sublinear expansion, and the many uses that have been found for it, we invite the reader to turn to the survey by Letzter~\cite{shohamsurvey} which also appears in this volume.
%Finding colour sets

While our discussion so far gives some indication of the approach taken in~\cite{montgomeryrbs} to study the algebraic properties of colourings, the actual implementation is somewhat more complicated. Aside from the omission here of many details around the distributive absorption approach, this is in part because the weight distribution on the auxiliary graph $K$ can vary widely. For example, when the colouring of $G$ arises directly from the addition table of an abelian group, then every weight will be 0. On the other hand, if the colouring of $G$ is a randomly chosen optimal colouring, then we expect $K$ to be roughly evenly weighted with total weight $\Theta(n^5)$, and, between these two examples, there are a wide variety of plausible auxiliary graphs $K$.

%While, if the $G$ arises from the `blow-up construction' of $\mathbb{Z}_2$ described in Section~\ref{sec:examples} and depicted in Figure~\ref{fig:moreeg}, with the colouring randomly chosen within the restriction given, then the auxiliary graph will (with high probability) consist of two cliques, on $C_0$ and $C_1$ respectively, where the edges are all roughly equally weighted.

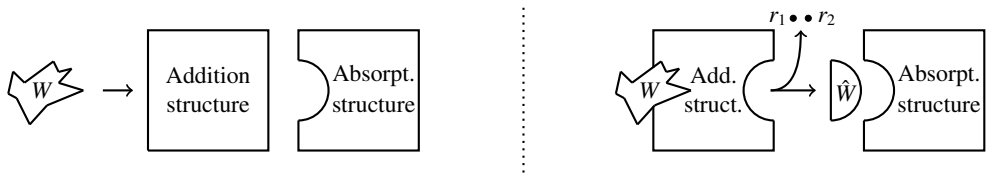
\begin{figure}
\begin{center}
\begin{tikzpicture}[scale=1]

\def\hgt{0.8};
\def\labh{0.2};

\draw [thick,->] ($(-0.2,0)-(\hgt,0)$) -- ($(0.2,0)-(\hgt,0)$);

\begin{scope}[shift={(0.5*\hgt,0)}]
\draw [white,thick,dotted] (0,1.1) -- (0,-1.1);

\draw [thick] (-\hgt,-\hgt) -- (-\hgt,\hgt) -- (\hgt,\hgt) -- (\hgt,-\hgt) -- (-\hgt,-\hgt);

\draw (0,\labh) node  {\footnotesize Addition};
\draw (0,-\labh) node {\footnotesize structure};
\end{scope}

\begin{scope}[shift={(3*\hgt,0)}]
\draw [thick] (-\hgt,-0.5*\hgt) -- (-\hgt,-\hgt) -- (\hgt,-\hgt) -- (\hgt,\hgt) -- (-\hgt,\hgt) -- (-\hgt,0.5*\hgt);

\draw [thick] (-\hgt,0.5*\hgt) arc(90:-90:0.5*\hgt) (-\hgt,-0.5*\hgt)  --cycle;

\draw ($(0.2,\labh)$) node  {\footnotesize Absorpt.};
\draw ($(0.2,-\labh)$) node {\footnotesize structure};
\end{scope}

\begin{scope}[shift={(-2.2*\hgt,0)}]

\coordinate (A1) at (0:0.5);
\coordinate (A2) at (30:0.2);
\coordinate (A3) at (40:0.45);
\coordinate (A4) at (55:0.3);
\coordinate (A5) at (75:0.4);
\coordinate (A5b) at (160:0.3);
\coordinate (A6) at (150:0.5);
\coordinate (A7) at (185:0.5);
\coordinate (A8) at (195:0.4);
\coordinate (A9) at (245:0.5);
\coordinate (A10) at (275:0.2);

\foreach \x/\y in {1/2,2/3,3/4,4/5,5/5b,5b/6,6/7,7/8,8/9,9/10,10/1}
{
\draw [thick] (A\x) -- (A\y);
}

\draw (-0.05,-0.0) node {\footnotesize $W$};
%\draw [thick] (A1) -- (A2) -- (A3) -- (A4) -- (A5) -- (A6)-- (A7) -- (A8) -- (A9) -- (A5) -- (A1);
\end{scope}

\end{tikzpicture}\hspace{1.2cm}\begin{tikzpicture}\draw [thick,dotted] (0,1.1) -- (0,-1.1);\end{tikzpicture}\hspace{1.2cm}\begin{tikzpicture}[scale=1]

\def\hgt{0.8};
\def\labh{0.2};

\draw [white,thick,dotted] (0,1.1) -- (0,-1.1);

%ADDDD
\draw [thick] (-\hgt,-0.21) -- (-\hgt,-\hgt) -- (\hgt,-\hgt) -- (\hgt,-0.5*\hgt);
\draw [thick] (\hgt,0.5*\hgt) -- (\hgt,\hgt) -- (-\hgt,\hgt) -- (-\hgt,0.2975);

\draw [thick] (\hgt,0.5*\hgt) arc(90:270:0.5*\hgt) (\hgt,-0.5*\hgt)  --cycle;
\draw (0,\labh) node  {\footnotesize Add.};
\draw (0,-\labh) node {\footnotesize struct.};

\begin{scope}[shift={(0.95*\hgt,0)}]
\draw [thick] (\hgt,0.5*\hgt) arc(90:-90:0.5*\hgt) (\hgt,-0.5*\hgt)  -- cycle;
\draw [thick]  (\hgt,0.5*\hgt) -- (\hgt,-0.5*\hgt);
\draw ($(\hgt+0.2,0.0)$) node {\footnotesize $\hat{W}$};
\draw [thick,->] (0,0) -- (0.75*\hgt,0);
\draw [thick,->] (0,0)  to[out=0,in=270] ($(0.75*\hgt-0.2,1*\hgt)$);
\draw ($(0.75*\hgt+0.15,1.2*\hgt)$) node {\footnotesize $r_2$};
\draw ($(0.75*\hgt-0.5,1.2*\hgt)$) node {\footnotesize $r_1$};
\draw [fill] ($(0.75*\hgt-0.1,1.2*\hgt)$) circle [radius=0.04cm];
\draw [fill] ($(0.75*\hgt-0.3,1.2*\hgt)$) circle [radius=0.04cm];

\end{scope}

\begin{scope}[shift={(3.5*\hgt,0)}]
\draw [thick] (-\hgt,-0.5*\hgt) -- (-\hgt,-\hgt) -- (\hgt,-\hgt) -- (\hgt,\hgt) -- (-\hgt,\hgt) -- (-\hgt,0.5*\hgt);

\draw [thick] (-\hgt,0.5*\hgt) arc(90:-90:0.5*\hgt) (-\hgt,-0.5*\hgt)  --cycle;

\draw ($(0.2,\labh)$) node  {\footnotesize Absorpt.};%\labh
\draw ($(0.2,-\labh)$) node {\footnotesize structure};
\end{scope}

\begin{scope}[shift={(-1*\hgt,0)}]

\coordinate (A1) at (0:0.5);
\coordinate (A2) at (30:0.2);
\coordinate (A3) at (40:0.45);
\coordinate (A4) at (55:0.3);
\coordinate (A5) at (75:0.4);
\coordinate (A5b) at (160:0.3);
\coordinate (A6) at (150:0.5);
\coordinate (A7) at (185:0.5);
\coordinate (A8) at (195:0.4);
\coordinate (A9) at (245:0.5);
\coordinate (A10) at (275:0.2);

\foreach \x/\y in {1/2,2/3,3/4,4/5,5/5b,5b/6,6/7,7/8,8/9,9/10,10/1}
{
\draw [thick] (A\x) -- (A\y);
}

\draw (-0.05,-0.0) node {\footnotesize $W$};
%\draw [thick] (A1) -- (A2) -- (A3) -- (A4) -- (A5) -- (A6)-- (A7) -- (A8) -- (A9) -- (A5) -- (A1);
\end{scope}

\end{tikzpicture}
\end{center}
\caption{An addition structure and an absorption structure on the left, along with a small unstructured vertex set $W$. On the right, the addition structure transforms the unstructured set into two `remainder vertices' and a structured vertex set $\hat{W}$ which is suitable for absorption.}\label{fig:add}
\end{figure}

Suppose, however, that our absorption structure can be created. Note that the condition on the set to be absorbed (that it is the vertex set of a matching of identity colour edges) is quite restrictive -- if a vertex is to be in the set to be absorbed then so must its neighbour across an identity-colour edge. To be able to consider an arbitrary balanced vertex set, an `addition structure' is introduced in~\cite{montgomeryrbs}. This consists of a set $V^{\mathrm{add}}$ of vertices and a set of colours $C^{\mathrm{add}}$ so that, given any sufficiently small set $W$ of vertices not in $V^{\mathrm{add}}$, $G$ contains two vertex-disjoint matchings $M_1$ and $M_2$ and two remainder vertices $r_1,r_2$, such that these vertices altogether exactly form the set $V^{\mathrm{add}}\cup W$, and such that $M_2$ is a rainbow matching using each colour in $C^{\mathrm{add}}$ exactly once and $M_1$ is a matching of identity colour edges. Thus, setting $\hat{W}=V(M_1)$, we have a set which can be absorbed by the absorption structure, while the matching $M_2$ is included in the rainbow matching that is found. The two remainder vertices are the two vertices not included in the final matching of $n-1$ edges. The function of the addition structure here, effectively, is to take a general unstructured set $W$ and transform it to a structured set $\hat{W}$ (with the loss of two remainder vertices) that is then suitable for absorption
(this is rather abstractly represented in Figure~\ref{fig:add}). The construction of the addition structure is not particularly complex, though we include no further details here, but we note that it works in an iterative fashion where more vertices are added to the structure in pairs, all while a changing pair of `remainder vertices' is considered along with two large matchings. Thinking of this as addition with remainder vertices is a very useful perspective when creating such a structure.

Finally, let us finish our discussion of the methods in~\cite{montgomeryrbs} by noting how useful it is to have two vertices which are omitted from the final matching. Considering the case where the colouring arises from an abelian group, given an arbitrary balanced set of vertices, in order for the addition structure to function any group element must be representable as the sum of two vertices from the addition structure (i.e., any remainder can be represented as the sum of two remainder vertices). This would be true if, say, the addition structure contained a random set of $n^{2/3}$ vertices, while using only one remainder vertex is not possible with these methods. Thus, perhaps it is fair to say that finding a rainbow matching in $G$ missing one edge is not as close  to an $n$-edge rainbow matching as it might seem. In some sense, the degree of the discrepancy is not 1 but 3: we may choose two vertices and one colour to leave out of the matching. For the methods in~\cite{montgomeryrbs} we need the flexibility to omit at least two of these --- a brief discussion of a variant of the problem which is approachable using this is in Section~\ref{sec:related:general}. Finding a rainbow matching without any of this flexibility (i.e., proving the Ryser-Brualdi-Stein conjecture for large odd $n$) seems quite far beyond these methods, and would certainly need new ideas.

 %%%%%%%%%%%%%%%%%%%%%%%%%%%%%%%%%%%%%%%%%%%%%%%%%%%%%%%%%%%%%%%%%%%%%%%%%%%%%%%%%%%%%%%%%%%%%%%%%%%%%%%%%%%%%%%%%%%%%%%%%%
 %%%%%%%%%%%%%%%%%%%%%%%%%%%%%%%%%%%%%%%%%%%%%%%%%%%%%%%%%%%%%%%%%%%%%%%%%%%%%%%%%%%%%%%%%%%%%%%%%%%%%%%%%%%%%%%%%%%%%%%%%%
 %%%%%%%%%%%%%%%%%%%%%%%%%%%%%%%%%%%%%%%%%%%%%%%%%%%%%%%%%%%%%%%%%%%%%%%%%%%%%%%%%%%%%%%%%%%%%%%%%%%%%%%%%%%%%%%%%%%%%%%%%%
 %%%%%%%%%%%%%%%%%%%%%%%%%%%%%%%%%%%%%%%%%%%%%%%%%%%%%%%%%%%%%%%%%%%%%%%%%%%%%%%%%%%%%%%%%%%%%%%%%%%%%%%%%%%%%%%%%%%%%%%%%%

\section{Related problems}\label{sec:related}%ainbow matchings in non-complete graphs and other r
In the rainbow subgraph formulation, we have discussed the problem of finding large rainbow matchings in optimally coloured complete bipartite balanced graphs. To finish this survey, we will discuss some related work and conjectures. We consider what happens when the colouring is proper but non-optimal (Section~\ref{sec:related:general}), when the coloured graph is not bipartite but complete (Section~\ref{sec:related:cycle}), when the colouring is not required to be proper (Section~\ref{sec:related:nonproper}) and
 when multigraphs are used (Section~\ref{sec:related:noncomplete}). Additionally, we discuss the closely related problem of large matchings in Steiner triple systems (Section~\ref{sec:related:STS}). We cover these topics beginning with those where the methods discussed in Section~\ref{sec:partial} give some progress.
Throughout this section, we touch only on the problems most closely connected to transversals in Latin squares, and, for additional background and a much wider variety of problems on rainbow subgraphs, we refer the reader to the excellent survey by Pokrovskiy~\cite{Alexeysurvey}.

\subsection{Generalised Latin squares}\label{sec:related:general}
From the perspective of the rainbow subgraph formulation, it is natural to ask how large a rainbow matching can be found in any properly coloured $K_{n,n}$. This corresponds to the study of transversals in \emph{generalised Latin squares}, also known as \emph{Latin arrays}, of order $n$, where an $n$ by $n$ grid is filled with symbols so that no symbol appears more than once in any row or column.  As non-optimal colourings have more colours than optimal colourings, it is generally expected that large, or even perfect, rainbow matchings should be easier to find in non-optimal colourings of $K_{n,n}$. For example, Montgomery, Pokrovskiy and Sudakov~\cite{montgomery2018decompositions} showed that any proper colouring of $K_{n,n}$ where at most $(1-o(1))n$ colours appear more than $(1-o(1))n$ times has a perfect rainbow matching.

The techniques used in work towards the Ryser-Brualdi-Stein conjecture have often also been applied more generally to proper colourings of $K_{n,n}$. For example, Keevash, Pokrovskiy, Sudakov and Yepremyan (see~\cite[Theorem 1.6]{KPSY}) showed that any such colouring contains a rainbow matching with $n-O(\log n/\log\log n)$ edges.
The results of~\cite{montgomery2018decompositions} show that, for large $n$, we only need consider colourings which are, in some sense, close to an optimal colouring. For sufficiently large $n$, this allows the techniques described in Section~\ref{sec:partial} to show that any properly coloured $K_{n,n}$ has a rainbow matching with $n-1$ edges~\cite{montgomeryrbs}.

Beyond this, Best, Pula and Wanless~\cite{best2021small} have conjectured that any properly coloured $K_{n,n}$ should have an $(n-1)$-edge rainbow matching despite the deletion of any chosen vertex, or, equivalently, that
every proper colouring of $K_{n-1,n}$ has a rainbow matching with $n-1$ edges.
As Georgakopoulos~\cite{georgakopoulos2013delay} observed, this is a special case of a conjecture of Haxell, Wilfong and Winkler (see also~\cite{alon2007edge}). Though not proved in full in~\cite{montgomeryrbs}, the methods described in Section~\ref{sec:partial} are strong enough to confirm this conjecture for sufficiently large $n$. (Very roughly, that such a partial transversal would omit at least one colour and omit one vertex is enough flexibility to allow the methods described in Section~\ref{sec:partial} to work, though it is clearer from the sketches in Section~\ref{sec:partial} that the methods would work in a setting where two vertices and no colours are omitted, such as to find an $(n-1)$-edge rainbow matching in an optimal colouring of $K_{n,n}$ with the edges of any one colour removed.) This conjecture is a weak version of a very strong conjecture by Stein~\cite{stein}, which we cover in Section~\ref{sec:related:noncomplete}.

%%%%%%%%%%%%%%%%%%%%%%%%%%%%%%%%%%%%%%%%%%%%%%%%%%%%%%%%%%%%%%%%%%%%%%%%%%%%%%%%%%%%%%%%%%%%%%%%

\subsection{Steiner triple systems}\label{sec:related:STS}
In addition to rainbow matchings in properly coloured balanced bipartite complete graphs, transversals in Latin squares also have a natural expression as a matching in 3-uniform hypergraphs. Given a Latin square $L$ of order $n$, form a hypergraph $\mathcal{H}(L)$ by first creating disjoint sets $A$ and $B$ of $n$ vertices representing the rows and columns of $L$ respectively and a set $C$ of $n$ vertices representing the symbols of $L$. Then, for each pair of vertices $a\in A$ and $b\in B$, take the symbol $c\in C$ in row $a$ and column $b$ of $L$ and add $abc$ to $\mathcal{H}(L)$.
Note that a partial transversal in the Latin square $L$ corresponds to a matching (a set of disjoint edges) in $\mathcal{H}(L)$, while a full transversal corresponds to a perfect matching.

Note further that any two vertices from different sets $A$, $B$ and $C$ are contained together in exactly one edge of $\mathcal{H}(L)$.
Thus, a Latin square of order $n$ can be represented as what is a tripartite version of a Steiner triple system.
A \emph{Steiner triple system (STS)} of order $n$ is a 3-uniform hypergraph with $n$ vertices in which every pair of vertices is contained in exactly one edge.
Steiner triple systems are a type of design whose history dates back to Kirkman's famous schoolgirl problem~\cite{kirkman} from 1850 (which asked for a particular type of Steiner triple system, called a \emph{resolvable design}, with 15 vertices). Every vertex in an STS of order $n$ must be in $\frac{n-1}{2}$ edges, and the number of vertex pairs, $\binom{n}{2}$, must be divisible by 3. Thus, Steiner triple systems of order $n$ can only exist if $n\equiv 1$ or $3\;\mathrm{mod}\; 6$. Kirkman~\cite{kirkman} showed that this condition is sufficient for the existence of an STS of order $n$.

%Originally, however, this was
%Matchings in hypergraphs more generally~\cite{rodlandhisnibble}.
The corresponding version of the Ryser-Brualdi-Stein conjecture was made for Steiner triple systems in 1981 by Brouwer~\cite{brouwer1981size}, as follows.

\begin{conj}[Brouwer]\label{conj:Brouwer}
Every Steiner triple system of order $n$ has a matching of at least $\frac{n-4}{3}$ edges.
\end{conj}

Conjecture~\ref{conj:Brouwer} would be tight for infinitely many $n$ as seen by constructions of Wilson (see~\cite{colbourn1992directed}) and Bryant and Horsley~\cite{bryant2013second,bryant2015steiner}.
Improving on early bounds by Wang \cite{wang} and Lindner and Phelps~\cite{lindner1978note},
 Brouwer~\cite{brouwer1981size} gave the first asymptotic version of Conjecture~\ref{conj:Brouwer} by showing that any STS of order $n$ has a matching with at least $n/3-O(n^{2/3})$ edges.
An asymptotic version of Conjecture~\ref{conj:Brouwer} can now, however, be proved with a simple application of standard matching theorems for almost-regular hypergraphs proved using the nibble method (see, for example,~\cite{alon2016probabilistic}). Refining such an approach and using large deviation inequalities, Alon, Kim and Spencer~\cite{AKS} improved Brouwer's result to show that a matching with $n/3-O(n^{1/2}\log^{3/2}n)$ edges exists in any STS of order $n$.

More recently, Keevash, Pokrovskiy, Sudakov and Yepremyan~\cite{KPSY} made the following conversion of this problem to one on rainbow matchings, using it to apply their methods and make a very significant improvement towards Conjecture~\ref{conj:Brouwer}. Assuming, for a slight simplication, that $n=3m$ for some integer $m$, take a Steiner triple system $S$ of order $n$ and partition $V(S)$ uniformly at random into sets $A,B$ and $C$, each with size $m$. Form a bipartite graph $G$ with vertex classes $A$ and $B$, and, for each $a\in A$ and $b\in B$, if there is some $c\in C$ such that $abc\in S$, then put an edge between $a$ and $b$ in $G$ with colour $c$.
 This gives a properly coloured bipartite graph, and, moreover, one in which a rainbow matching corresponds to a matching in the original Steiner triple system. The graph $G$ is not a complete bipartite graph, but does look roughly like an optimal colouring of $K_{n,n}$ where edges have been deleted independently at random with probability $2/3$. This allowed Keevash, Pokrovskiy, Sudakov and Yepremyan~\cite{KPSY} to hugely improve the previous bound to show that any STS of order $n$ has a matching with
$n/3-O(\log n/\log\log n)$ edges. Applying the techniques described in Section~\ref{sec:partial} through this conversion introduces some complications, but in~\cite{montgomeryrbs} it is shown that Conjecture~\ref{conj:Brouwer} holds for sufficiently large $n$.

As mentioned in Section~\ref{sec:decomp}, we can also consider perfect matchings in random Steiner triple systems of order $n\equiv 3\;\mathrm{mod}\; 6$. Kwan~\cite{kwan2020almost} showed that such a random STS of order $n\equiv 3\;\mathrm{mod}\; 6$ does have a perfect matching with high probability, and, moreover contains $\left((1-o(1))\frac{n}{2e^2}\right)^{n/3}$ perfect matchings with high probability. Morris~\cite{morris2017random} adapted Kwan's proof to show that a random STS of order $n\equiv 3\;\mathrm{mod}\; 6$ has $\Omega(n)$ disjoint perfect matchings with high probability. Ferber and Kwan~\cite{ferbkwan} were then able to show (see Section~\ref{sec:decomp}) that $(1-o(1))n$ disjoint perfect matchings can be found with high probability, while making the corresponding conjecture to Conjecture~\ref{conj:random}, as follows.

\begin{conj}[Ferber and Kwan]\label{conj:FK} A random Steiner triple system of order $n\equiv 3\;\mathrm{mod}\; 6$ has a decomposition into perfect matchings with high probability.
\end{conj}

A Steiner triple system with such a decomposition as in Conjecture~\ref{conj:FK} is known as \emph{resolvable}. Whether there existed any resolvable STSs for each $n\equiv 3\;\mathrm{mod}\; 6$ was an old problem in combinatorics. Their existence was known for several infinite families of integers $n$ (see~\cite{hall1967combinatorial}), but in general was only shown  in famous work by
Ray-Chaudhuri and Wilson~\cite{ray1971solution} in 1971.

%%%%%%%%%%%%%%%%%%%%%%%%%%%%%%%%%%%%%%%%%%%%%%%%%%%%%%%%%%%%%%%%%%%%%%%%%%%%%%%%%%%%%%%%%%%%%%%%

\subsection{Complete graphs and rainbow {H}amilton cycles}\label{sec:related:cycle}
What happens if we, instead, look for a large rainbow matching in a properly coloured $n$-vertex complete graph, $K_n$?  The limit on the size of the matching imposed by the number of vertices ($\lfloor n/2\rfloor$) is much stronger than that imposed by the number of colours needed in a proper colouring ($n-1$ when $n$ is even). This makes the problem much easier, and it is more interesting
to ask for a rainbow subgraph with more edges, or even $n-1$ edges, so that almost all the colours might be used. The most natural subgraph to consider is a Hamilton path. (For the consideration of trees more generally, we refer the reader to the survey by Pokrovskiy~\cite{Alexeysurvey}.) However, giving a counterexample to a conjecture by Hahn~\cite{hahn1980jeu}, in 1984 Maamoun and Meyniel~\cite{maamoun1984problem} gave an optimal colouring of $K_n$ which has no rainbow Hamilton path when $n\geq 4$ is any power of $2$. To see this colouring, take the complete graph with the group $\mathbb{Z}_2^k$ (with $k\geq 2$) as its vertex set, where each edge $xy$ is coloured $x+y$. Using that $2x=0$ for every $x\in \mathbb{Z}_2^k$,
this can easily be seen to be a proper colouring which never uses 0, so that the colouring with $\mathbb{Z}_2^k\setminus\{0\}$
is optimal. If this graph has a rainbow Hamilton path, then the corresponding equation to \eqref{eqn:sum} shows that the sum of the two endvertices must be 0, a contradiction. Therefore, the following conjecture of Andersen~\cite{andersen1989hamilton} from 1989 would be best possible when $n\geq 4$ is a power of 2.

\begin{conj}[Andersen] All proper edge-colourings of $K_n$ have a rainbow path with length $n-2$.\label{conj:and}
\end{conj}

This question can also be asked for rainbow cycles, where we note that, when $n$ is even, any optimal colouring of $K_n$ does not have a rainbow Hamilton cycle simply as it does not have enough colours to support this.
Akbari, Etesame, Mahini and Mahmoody~\cite{akbari2007rainbow} have asked whether all optimal colourings of $K_n$ have a rainbow cycle with length $n-2$, although the current author knows of no optimal colouring of $K_n$ with no rainbow cycle with length $n-1$. We also can note that the known optimal colourings of $K_n$ with no rainbow Hamilton path are much rarer than the many examples of Latin squares with no full transversals given in Section~\ref{sec:examples} (and are only known for when $n\geq 4$ is a power of 2). It is, however, known that a random optimal colouring of $K_n$ will have a rainbow Hamilton path with high probability, due to Gould, Kelly, K\"uhn and Osthus~\cite{gould2022almost}, who also showed that such a random colouring will have a rainbow cycle containing all of the colours with high probability.

To date, the progress towards Conjecture~\ref{conj:and} has not needed to distinguish between the search for long rainbow paths or cycles. In contrast to rainbow matchings in the bipartite case, it turns out to be difficult even to find a rainbow path or cycle containing $(1-o(1))n$ vertices, and (after various improved bounds, as detailed in~\cite{alon2017random}) this was only done in 2017, by
Alon, Pokrovskiy and Sudakov~\cite{alon2017random}. More precisely, in~\cite{alon2017random} it is shown that any proper colouring of $K_n$ contains a rainbow cycle with $n-O(n^{3/4})$ vertices. Balogh and Molla~\cite{balogh2019long} have used the techniques in~\cite{alon2017random} to show that in fact such a cycle with $n-O(\sqrt{n}\log n)$ vertices always exists.

Techniques used to work towards the Ryser-Brualdi-Stein conjecture do not seem to translate naturally to work on Conjecture~\ref{conj:and}, and, for example, Keevash, Pokrovskiy, Sudakov and Yepremyan~\cite{KPSY} did not apply their techniques to this problem (see~\cite{Alexeysurvey}). Similarly, the techniques described in Section~\ref{sec:partial} do not apply simply to work towards Conjecture~\ref{conj:and}. However, with significant development, it seems likely progress can be made by considering not only a single special `identity colour', but a basis of colours that can be used to represent the full colour set, while carrying out a more complicated addition using these representations. This seems likely to show that, in any proper colouring of $K_n$, there is a rainbow cycle with $n-O(1)$ vertices and is the subject of forthcoming work by Benford, Bowtell and Montgomery~\cite{BBM}.

\subsection{Non-proper colourings}\label{sec:related:nonproper}
Interestingly, that the colourings we have considered are proper may be an unnecessarily strong condition if we want only to find a very large rainbow matching. In particular, in 1975 Stein~\cite{stein} made a series of 7 conjectures on transversals in arrays where the restrictions of Latin arrays are replaced (at least in part) by a restriction only on the number of times a symbol appears in the entire array (or correspondingly only a restriction on the number of times a colour appears in our edge-coloured graphs).
Granted, it turns out that these conjectures are mostly false, as shown by a counterexample by Pokrovskiy and Sudakov~\cite{pokrovskiy2019counterexample} and small variants of this (see~\cite{best2021small}), but there remains the potential that very large partial transversals can still be found under the same conditions.

For example, Stein conjectured that every equi-$n$ square should have a partial transversal of order $n-1$, where an \emph{equi}-$n$ square is an $n$ by $n$ square in which each one of $n$ symbols appears exactly $n$ times.
Pokrovskiy and Sudakov~\cite{pokrovskiy2019counterexample} constructed an equi-$n$ square with no partial transversal of order $n-\frac{1}{42}\log n$, but it remains an open question whether such squares contain a partial transversal of order $(1-o(1))n$. In this direction, the best bound known is by Aharoni, Berger, Kotlar, and
Ziv~\cite{aharoni2017conjecture}, who used topological methods to show that every equi-$n$-square has a partial transversal of order at least
$2n/3$.

Finally, here, let us note that one of Stein's conjectures in~\cite{stein} does remain open (see~\cite{best2021small}). In the language of coloured graphs, this states that any colouring of $K_{n,n-1}$ in which no vertex in the larger vertex class is adjacent to more than one edge of any colour (i.e., the colouring is proper on one side) should contain a rainbow matching with $n-1$ edges. This conjecture is much stronger than the special case conjectured by Best, Pula and Wanless~\cite{best2021small} mentioned in Section~\ref{sec:related:general}, and, as such, remains widely open.

\subsection{Multigraphs}\label{sec:related:noncomplete}
%Finally, let us turn to properly coloured graphs which are neither complete nor complete bipartite.
%Answering a problem posed by Wang~\cite{wang2011rainbow}, Diemunsch, Ferrara, Moffatt, Pfender, and Wenger~\cite{diemunsch2011rainbow} proved that we can define some function $f(d)$ such that any properly coloured graph with minimum degree at least $f(d)$ has a rainbow matching of $d$ edges (and, indeed, that $f(d)\leq $).
%Lo and Tan~\cite{lo2014note}
%Gy\'arf\'as and Sark\"ozy~\cite{gyarfas2012rainbow}
Finally, let us turn to an interesting generalisation of the Ryser-Brualdi-Stein conjecture. Considering an optimal colouring of $K_{n,n}$, we have $n$ monochromatic (disjoint) matchings of size $n$, and wish to form a matching with $n$ or $n-1$ edges using at most one edge per matching. What if we no longer require these matchings to be disjoint, so that together they form a multigraph? This, in an even more general form, is considered by the following conjecture of Aharoni and Berger~\cite{aharoni2009rainbow} (where $n$ corresponds to $n-1$ in our discussion so far).

\begin{conj}
\label{conj:AB}
If $G$ is a bipartite multigraph formed of $n$ matchings of size $n+1$, each of a different colour, then $G$ contains a rainbow matching with $n$ edges.
\end{conj}

After a sequence of improvements by various authors (see, for example,~\cite{Alexeysurvey}), the best currently known bound towards Conjecture~\ref{conj:AB} that holds for all $n$ is by Aharoni, Kotlar,
and Ziv~\cite{aharoni2017representation}, who show that if the matchings in Conjecture~\ref{conj:AB} are made large, with at least $3n/2+1$ edges each, then there is always an $n$-edge rainbow matching. Conjecture~\ref{conj:AB} has also been proven asymptotically by Pokrovskiy~\cite{pokrovskiy2018approximate}, who showed that if there are $(1+o(1))n$ edges in each matching then there is an $n$-edge rainbow matching. Though this asymptotic result now has a considerably easier proof, by Munh\'a Correia, Pokrovskiy and Sudakov~\cite{correia2023short}, the full conjecture remains ambitious. The techniques described in Section~\ref{sec:partial} do not seem to be applicable and, thus, Conjecture~\ref{conj:AB} remains very open indeed.

\thankyou{The author would like to thank colleagues, in particular Ian Wanless and the anonymous referee, for their comments on this survey. The author would also like to acknowledge support by the European Research Council (ERC) under the European Union Horizon 2020 research and innovation programme (grant agreement No. 947978) and the Leverhulme trust. }%Put any acknowledgements here.}

%You can also use BibTeX with the amsplain style.
\bibliographystyle{abbrv}

%\bigskip
%For preference, use the BibTeX style \texttt{amsplain} and send the editors
%your \texttt{.bib} file with exactly the references you use.

\myaddress

%% Editors to uncomment the two lines below if survey ends on an odd-numbered page
%\newpage
%\mbox{}

\end{document}